\documentclass[10pt,twoside]{article}
\usepackage[utf8]{inputenc}

\title{\uppercase{\textbf{\large{Beilinson--Lichtenbaum phenomenon for motivic cohomology}}}}

\author{TESS BOUIS AND ARNAB KUNDU}
\date{}
\usepackage[left=2.25cm, right=2.25cm, top=1in, bottom=1in, headheight=2cm]{geometry}

\usepackage{fancyhdr}
\usepackage{stmaryrd}

\usepackage{mathdots}

\usepackage{mathtools}

\usepackage{amsmath, amscd, amsfonts, amssymb, amsthm,todonotes,eurosym,enumitem,relsize}

\usepackage[utf8]{inputenc}
\usepackage{dirtytalk}
\usepackage[T1]{fontenc}
\usepackage{mathrsfs}

\usepackage[colorlinks]{hyperref}
\usepackage[figure,table]{hypcap}
\definecolor{imperialred}{RGB}{237, 41, 57}
\definecolor{royalblue}{RGB}{64, 106, 212}
\definecolor{link}{RGB}{11,0,128}
\definecolor{gren}{RGB}{32,130,63}
\hypersetup{
	bookmarksnumbered,
	pdfstartview={FitH},
    citecolor=royalblue,
 	linkcolor=imperialred,
 	urlcolor=link,
 	linktocpage = true
	pdfpagemode={UseOutlines}
}

\usepackage{xcolor}
\usepackage{comment}

\usepackage[object=pgfhan]{pgfornament}

\usepackage{graphicx}
\usepackage{array}
\usepackage{csquotes}
\usepackage{extarrows}

\usepackage{needspace}

\input{xy}
\xyoption{all}
\usepackage{tikz-cd}
\usepackage{textcomp}
\usepackage{colonequals}

\DeclareFontFamily{T1}{cbgreek}{}
\DeclareFontShape{T1}{cbgreek}{m}{n}{<-6>  grmn0500 <6-7> grmn0600 <7-8> grmn0700 <8-9> grmn0800 <9-10> grmn0900 <10-12> grmn1000 <12-17> grmn1200 <17-> grmn1728}{}
\DeclareSymbolFont{quadratics}{T1}{cbgreek}{m}{n}
\DeclareMathSymbol{\qoppa}{\mathord}{quadratics}{19}
\DeclareMathSymbol{\Qoppa}{\mathord}{quadratics}{21}

\DeclareMathAlphabet\matheu{U}{eus}{m}{n}
\DeclareMathAlphabet{\mathmm}{U}{mmambb}{m}{n}
\newcommand{\bb}[1]{\mathbb{#1}}

\newcommand{\ca}[1]{\mathcal{#1}}
\newcommand{\scr}[1]{\mathscr{#1}}

\newcommand{\bun}{\begin{itemize}}
\newcommand{\bn}{\begin{enumerate}}
\newcommand{\bt}{\begin{thm}}
\newcommand{\bl}{\begin{lem}}
\newcommand{\bp}{\begin{prop}}
\newcommand{\bc}{\begin{cor}}
\newcommand{\bd}{\begin{teqn}}
\newcommand{\bud}{\begin{teqn*}}	
\newcommand{\bs}{\begin{proof}}
\newcommand{\br}{\begin{rem}}
\newcommand{\bdf}{\begin{defn}}
\newcommand{\bcj}{\begin{conj}}

\newcommand{\eun}{\end{itemize}}
\newcommand{\en}{\end{enumerate}}
\newcommand{\et}{\end{thm}}
\newcommand{\el}{\end{lem}}
\newcommand{\ep}{\end{prop}}
\newcommand{\ec}{\end{cor}}
\newcommand{\ed}{\end{teqn}}
\newcommand{\eud}{\end{teqn*}}
\newcommand{\es}{\end{proof}}
\newcommand{\er}{\end{rem}}
\newcommand{\edf}{\end{defn}}
\newcommand{\ecj}{\end{conj}}

\DeclareMathOperator{\codim}{codim}

\DeclareMathOperator{\Frac}{Frac}
\DeclareMathOperator{\spec}{Spec}



\usepackage{sectsty}

\usepackage{fancyhdr}

\newlength{\outermargin} 
\newlength{\mar} 
\newlength{\len}
\newlength{\temp}
\newtheorem{theorem}{Theorem}[section]
\newtheorem{lemma}[theorem]{Lemma}
\newtheorem{corollary}[theorem]{Corollary}
\newtheorem{proposition}[theorem]{Proposition}

\newtheorem{theoremintro}{Theorem}

\newtheorem{corollaryintro}[theoremintro]{Corollary}

\theoremstyle{definition}
\newtheorem{remark}[theorem]{Remark}

\newtheorem{examples}[theorem]{Examples}

\newtheorem{definition}[theorem]{Definition}

\DeclareMathOperator{\Z}{\mathbb{Z}}
\DeclareMathOperator{\Q}{\mathbb{Q}}

\DeclareMathOperator{\N}{\mathbb{N}}
\DeclareMathOperator{\F}{\mathbb{F}}

\renewcommand{\epsilon}{\varepsilon}

\DeclareFontFamily{U}{MnSymbolC}{}
\DeclareFontShape{U}{MnSymbolC}{m}{n}{
	<-5.5> MnSymbolC5
	<5.5-6.5> MnSymbolC6
	<6.5-7.5> MnSymbolC7
	<7.5-8.5> MnSymbolC8
	<8.5-9.5> MnSymbolC9
	<9.5-11.5> MnSymbolC10
	<11.5-> MnSymbolCb12
}{}

\usepackage[bbgreekl]{mathbbol}
\DeclareSymbolFontAlphabet{\mathbb}{AMSb}
\DeclareSymbolFontAlphabet{\mathbbl}{bbold}

\newcommand{\ccite}[2]{\cite[#2]{#1}}
\newcommand{\stacks}[1]{\cite[\href{https://stacks.math.columbia.edu/tag/#1}{Tag #1}]{stacks-project}}
\usepackage{soul}
\usepackage{enumitem}
\numberwithin{equation}{theorem}

\begin{document}

\maketitle
	
	\pagestyle{fancy}
	\fancyhead[EC]{TESS BOUIS AND ARNAB KUNDU}
	\fancyhead[OC]{\uppercase{Beilinson--Lichtenbaum phenomenon for motivic cohomology}}
	\fancyfoot[C]{\thepage}
	
\begin{abstract}
The goal of this paper is to study non-$\mathbb{A}^1$-invariant motivic cohomology, recently defined by Elmanto, Morrow, and the first-named author, for smooth schemes over possibly non-discrete valuation rings. We establish that the cycle class map from $p$-adic motivic cohomology to a suitable truncation of Bhatt--Lurie's syntomic cohomology is an isomorphism, thereby verifying the Beilinson--Lichtenbaum conjecture in this generality. As a first consequence, we prove that this motivic cohomology integrally recovers the classical definition of motivic cohomology in terms of Bloch's cycle complexes, whenever the latter is defined. As a second consequence, we show a purity theorem for this cohomology theory over perfectoid rings, thus motivically refining a result of Nizio\l{} in algebraic $K$-theory. The key ingredient in our approach is a version of Gabber's presentation lemma applicable in mixed characteristic, non-noetherian settings.
\end{abstract}

{
  \hypersetup{linkcolor=black}
  \tableofcontents
}

\section{Introduction}

\vspace{-\parindent}
\hspace{\parindent}

Beilinson and Lichtenbaum, inspired by Grothendieck's vision of motives, conjectured the existence of a theory of {\it motivic cohomology} satisfying a list of properties, including the expression of this cohomology theory as the graded pieces of a filtration on algebraic 
$K$-theory \cite{lichtenbaum_values_1973,lichtenbaum_values_1984,beilinson_notes_1986,beilinson_height_1987,beilinson_notes_1987}.

Using the framework of motivic $\mathbb{A}^1$-homotopy theory, which ultimately led to the proof of the Bloch--Kato conjecture \cite{voevodsky_motivic_2011}, Voevodsky initiated the development of this cohomology theory for smooth schemes over a field \cite{bloch_algebraic_1986,voevodsky_cycles_2000}. This theory was later generalised to smooth schemes over any Dedekind domain by a number of authors \cite{levine_techniques_2001,geisser_motivic_2004,spitzweck_commutative_2018,bachmann_very_2025}. In both situations, the authors showed that {\it Bloch's cycle complexes} $z^i(X,\bullet-2i)$ satisfy most of Beilinson and Lichtenbaum's conjectures. 

Recently, Elmanto and Morrow \cite{elmanto_motivic_2023}, using trace methods in algebraic $K$-theory, introduced a new definition of motivic cohomology that applies to arbitrary quasi-compact quasi-separated (qcqs) schemes over a field. Building on their work, the first-named author of the present article extended this approach to define motivic complexes $\Z(i)^{\text{mot}}(X)$ for arbitrary qcqs schemes $X$ over $\bb{Z}$ \cite{bouis_motivic_2024,bouis_weibel_2025}. This theory, which satisfies various expected structural properties of motivic cohomology, may thus provide a more flexible framework for studying algebraic cycles. To make this expectation precise, it is therefore natural to ask how this theory compares to the more classical theory of Bloch's cycle complexes in the smooth case. 

Our first result, given below, is to show that the complexes $\Z(i)^{\text{mot}}$ coincide with Bloch's cycle complexes in the smooth case, {\it i.e.}, that the motivic complexes $\Z(i)^{\text{mot}}$ extend the classical definition of motivic cohomology beyond the case of smooth schemes. In particular, this suggests that the motivic complexes $\Z(i)^{\text{mot}}$ are well-suited to understand the interplay between algebraic cycles in characteristic zero and in positive characteristic.

\begin{theoremintro}[See Theorem~\ref{theoremclassicalmotiviccomparison}]\label{theoremintrocomparisonclassicalmotivic}
    Let $X$ be an ind-smooth scheme over a Dedekind domain. 
    Then for every integer $i \geqslant 0$, the natural comparison map
        $$z^i(X,\bullet-2i) \longrightarrow \Z(i)^{\emph{mot}}(X)$$
    is an equivalence in the derived category $\mathcal{D}(\Z)$ after Zariski sheafification.
\end{theoremintro}

For smooth schemes over a field, this result was proved by Elmanto--Morrow \cite{elmanto_motivic_2023}, using Gabber's presentation lemma \cite{colliot-thelene_bloch_1997}. Our proof of Theorem~\ref{theoremintrocomparisonclassicalmotivic} uses a Nisnevich-local variant of Gabber's presentation lemma in mixed characteristic (Theorem~\ref{theoremintropresentationlemma} below) to reduce the statement to the case where $X$ itself is the spectrum of a discrete valuation ring. In this case, the equivalence follows from a special case of the Beilinson--Lichtenbaum conjecture proved by Geisser \cite{geisser_motivic_2004} as a consequence of the Bloch--Kato conjecture, together with a special case of Theorem~\ref{theoremintroBL} established in \cite{bouis_motivic_2024}, by comparing both sides to the syntomic cohomology of $X$.

Our next results concern possibly non-noetherian rings. These rings---most prominently perfectoid rings \cite{scholze_perfectoid_2012,bhatt_arc-topology_2021,bhatt_prisms_2022}---have become a powerful tool in modern $p$-adic geometry. In many contemporary approaches, results for noetherian objects are indeed proved by first reducing claims to suitable non-noetherian ones. This shift provides the main motivation for studying algebraic geometry over schemes that are not necessarily noetherian.

We fix a prime number $p$ for the rest of this introduction. We also let $L_{\text{Nis}}$ denote the Nisnevich sheafification functor, and $\Z/p^k(i)^{\text{syn}}$ denote the syntomic cohomology complexes, as defined by Bhatt and Lurie \cite{bhatt_absolute_2022}. By construction, these syntomic complexes provide a generalisation of the $p$-adic étale cohomology complexes when $p$ is not necessarily invertible in the input.

The following result generalises (and relies upon) the Beilinson--Lichtenbaum comparison, stating that motivic cohomology can be expressed as a suitable truncation of syntomic cohomology, from the classical setting where the base is a field or a discrete valuation ring to that of an arbitrary valuation ring.\footnote{A valuation ring is an integral domain $V$ such that for any elements $f$ and $g$ in $V$, either $f \in gV$ or $g \in fV$. Note that a valuation ring is noetherian if and only if it is a discrete valuation ring.} Note that this comparison is not expected to hold over an arbitrary base ring (see Remark~\ref{remarkBLfalseingeneral}).

\setcounter{footnote}{0}
\begin{theoremintro}[Non-noetherian Beilinson--Lichtenbaum conjecture; see Theorem~\ref{theoremBLPrüfer}]\label{theoremintroBL}
    Let $X$ be an ind-smooth scheme over a valuation ring. 
    Then for all integers $i \geqslant 0$ and $k \geqslant 1$, there is a natural equivalence
    $$\Z/p^k(i)^{\emph{mot}}(X) \simeq \big(L_{\emph{Nis}} \tau^{\leqslant i} \Z/p^k(i)^{\emph{syn}}\big)(X)$$
    in the derived category $\mathcal{D}(\Z/p^k)$.
\end{theoremintro}

As a consequence of Theorem~\ref{theoremintrocomparisonclassicalmotivic}, for smooth schemes over discrete valuation rings, Theorem~\ref{theoremintroBL} agrees with the classical statement of the Beilinson--Lichtenbaum conjecture, proved in \cite{geisser_motivic_2004,voevodsky_motivic_2011}.

The proof of the Beilinson--Lichtenbaum comparison for smooth schemes over fields in \cite{voevodsky_motivic_2011} relies (among many things) on Gersten resolutions in order to reduce the statement to the case of fields (see \cite{geisser_bloch-kato_2001}). Our proof of Theorem~\ref{theoremintroBL} relies on a Nisnevich-local variant of Gabber's presentation lemma which holds over an arbitrary valuation ring (Theorem~\ref{theoremintropresentationlemma} below), which we use to reduce the statement to the case of henselian valuation rings. In this case, the result is a consequence of \cite[Proposition~$5.3$]{bouis_motivic_2026}.

As a consequence of Theorem~\ref{theoremintroBL}, we obtain the following description of motivic cohomology of smooth schemes over valuation rings in terms of the $p$-adic étale cohomology of their generic fibre.

\par To formulate this result, we use the notion of a $F$\nobreakdash-smooth ring due to Bhatt and Mathew \cite{bhatt_syntomic_2023}, which is a non-noetherian generalisation of the notion of a regular ring. Examples of $F$-smooth rings include regular rings (in particular, discrete valuation rings), perfectoid rings, and any valuation ring containing a field or whose $p$-completion contains a perfectoid valuation ring \cite{bhatt_syntomic_2023,kelly_k-theory_2021,bouis_cartier_2023}. Conjecturally, all valuation rings should be $F$-smooth (see Remark~\ref{remarkconjecturevaluationringsFsmooth}).

\begin{theoremintro}\label{theoremintroétale}
    Let $V$ be a $p$-torsionfree $F$-smooth valuation ring ({\it e.g.}, the ring of integers of an algebraic extension of a $p$-adic local field), and $R$ be a henselian local ind-smooth algebra over $V$. Then, for all integers $i,n \geqslant 0$ and $k \geqslant 1$, there is a natural isomorphism
    $$\emph{H}^n_{\emph{mot}}(\emph{Spec}(R),\Z/p^k(i)) \cong \begin{cases} 
        \emph{H}^n_{\emph{ét}}(\emph{Spec}(R[\tfrac{1}{p}]),\mu_{p^k}^{\otimes i}) 
    & \text{if } n < i \\
        \emph{Im}((R^\times)^{\otimes i} \longrightarrow \emph{H}^n_{\emph{ét}}(\emph{Spec}(R[\tfrac{1}{p}]),\mu_{p^k}^{\otimes i}) )    & \text{if } n = i\\
0 
    & \text{if } n > i  
\end{cases}$$
    of abelian groups.
\end{theoremintro}

Our proof of Theorem~\ref{theoremintroétale} proceeds in two steps: first, we use Theorem~\ref{theoremintroBL} to compare motivic cohomology with a truncation of syntomic cohomology; we then apply results of Bhatt and Mathew \cite{bhatt_syntomic_2023} in integral $p$-adic Hodge theory to compute this truncated syntomic cohomology in terms of étale cohomology. 
\par As a consequence of Theorem~\ref{theoremintroétale}, we obtain the following corollary, which establishes a form of purity for motivic cohomology over a perfectoid base. We note that, in particular, our proof of Corollary~\ref{corollaryintrocomparisongenericfibre} does not rely on absolute purity (see Remark~\ref{rem:annala_elmanto}).

\begin{corollaryintro}[See Corollary \ref{corollarycomparisongenericfibre}]\label{corollaryintrocomparisongenericfibre}
    Let $V$ be a perfectoid valuation ring of residue characteristic $p$, $F$ be the fraction field of $V$, and $X$ be an ind-smooth scheme over $V$. Then for all integers $i \geqslant 0$ and $k \geqslant 1$, the natural map
    $$\Z/p^k(i)^{\emph{mot}}(X) \longrightarrow \Z/p^k(i)^{\emph{mot}}(X_F)$$
    is an equivalence in the derived category $\mathcal{D}(\Z/p^k)$.
\end{corollaryintro}

The case of the above result where $V = \mathcal{O}_{\overline{K}}$, for $K$ a $p$-adic local field, is a cohomological refinement of a $K$-theoretic result proved by Nizio\l{} to prove Fontaine's crystalline conjecture \cite{niziol_crystalline_1998}. This cohomological result can in turn be used to give a more streamlined version of this proof (see Remark~\ref{remarkcrystallineconjecture}).

As a final application of Theorem~\ref{theoremintroBL}, we demonstrate that the motivic complexes $\Z(i)^{\text{mot}}$ are typically $\mathbb{A}^1$\nobreakdash-invariant on smooth schemes over valuation rings (see Theorem~\ref{theoremA1invariancemain}), using recent results of Bachmann--Elmanto--Morrow \cite{bachmann_A^1-invariant_2025}. This partially answers a question of Antieau--Mathew--Morrow \cite{antieau_K-theory_2022}.

We end this introduction with the following version of Gabber's presentation lemma over a general base, which serves as a technical basis for our Theorems~\ref{theoremintrocomparisonclassicalmotivic} and~\ref{theoremintroBL}.

\begin{theoremintro}[Presentation lemma;  see Corollary~\ref{cor:presentation_lemma}]\label{theoremintropresentationlemma}
    Let $S$ be the spectrum of a henselian local ring, $X$ be a smooth $S$\nobreakdash-scheme, and $Z \hookrightarrow X$ be a closed immersion of $S$-schemes which is $S$-fibrewise of positive codimension. Then for every point $x \in X$, there exist a Nisnevich neighbourhood $x\in X' \rightarrow X$, a smooth $S$-scheme $T$, and a Nisnevich square of $T$-schemes
    $$\begin{tikzcd}
        X' \setminus Z' \arrow[d]\arrow[r, hook] & X' \arrow[d, "\psi"]\\ 
        \mathbb{A}^1_T \setminus \psi(Z') \arrow[r, hook] & \mathbb{A}^1_T
    \end{tikzcd}$$
    where $Z' \colonequals Z \times_X X'$, $\psi$ is étale, and such that $\psi(Z')$ is finite over $T$.
\end{theoremintro}

Note that Gabber's presentation lemma over an infinite field \cite{gabber_gersten_conjecture} has already been extended to arbitrary fields by Hogadi--Kulkarni \cite{hogadi_kulkarni_gabber_finite_field}, to Dedekind domains with infinite residue fields by Schmidt--Strunk \cite{schmidt_strunk}, and to noetherian domains with infinite residue fields by Deshmukh--Hogadi--Kulkarni--Yadav \cite{deshmukh-hogadi-kulkarni-yadav_presentation_lemma_noetherian_domains}. A weaker form of Theorem~\ref{theoremintropresentationlemma}, where the finiteness of $\psi(Z')$ over $T$ is not guaranteed, was also obtained by Druzhinin \cite[Section~$4.2$]{druzhinin_a1_connectivity}. 

Our proof of Theorem~\ref{theoremintropresentationlemma} substantially simplifies the proofs of the presentation lemmas given in \cite{schmidt_strunk,deshmukh-hogadi-kulkarni-yadav_presentation_lemma_noetherian_domains,druzhinin_a1_connectivity}. In particular, our argument works over a base which is not necessarily noetherian, with no hypothesis on its residue fields, and we do not rely on the subtle intermediary construction of good compactifications for closed immersions. Finally, note that Theorem~\ref{theoremintropresentationlemma} is used in a crucial way in the sequel \cite{bouis_A^1-connectivity_2026} to this article to prove a version of Morel's unstable $\mathbb{A}^1$-connectivity theorem over arbitrary base schemes. 

\subsection*{Notation.}

\vspace{-\parindent}
\hspace{\parindent}

We will denote by $\text{Set}_{\ast}$ the $1$-category of pointed sets, by $[\ast]$ the trivial class of a pointed set, and by $\text{Sp}$ the stable $\infty$-category of spectra. We always index complexes using the cohomological grading.

A presheaf is called {\it finitary} if it commutes with filtered colimits of rings. For the Grothendieck topologies~$\tau$ that we will use, we denote by $L_{\tau}$ the associated sheafification functors. Given a scheme $X$ and a $\mathcal{D}(\Z)$-valued finitary $\tau$-sheaf $F$ on the small $\tau$-site of $X$, we say that $F$ lives $\tau$-locally on $X$ in degrees at most $i$ if $F$ takes values in $\mathcal{D}^{\leqslant i}(\Z)$ on the $\tau$-local rings of $X$. Given a point $x$ in a scheme $X$, by a \textit{Zariski neighbourhood} (resp. a \textit{Nisnevich neighbourhood}) of $x \in X$, we mean an open subset $x\in X' \subseteq X$ (resp. an étale morphism $X' \to X$ together with a lift of the point $x$ in $X'$). In either case, when $X'$ is moreover affine, we call this neighbourhood an affine neighbourhood. 

Unless otherwise specified, by dimension of a scheme, we mean its Krull dimension. A morphism of schemes $f\colon X\to S$ is said to be \textit{fibrewise of dimension $d$} (resp. \textit{fibrewise of dimension $\leqslant d$}), for some integer $d\geqslant 0$, if every fibre of $f$ is of dimension $d$ (resp. of dimension $\leqslant d$). The dimension of $f$ is defined to be the dimension of its fibre with the maximum dimension. 

\subsection*{Acknowledgements.}

\vspace{-\parindent}
\hspace{\parindent}

Thanks to the organisers of the Winter Workshop on $p$-adic Geometry at IAS, Princeton where this project started. The authors would like to express their sincere gratitude to Elden Elmanto and Matthew Morrow for their constant support and for sharing many insights. We would also like to thank Kęstutis Česnavičius, Frédéric Déglise, Neeraj Deshmukh, Marc Levine, Morten Lüders, Mura Yakerson, and Alex Youcis for helpful discussions, as well as Neeraj Deshmukh, Elden Elmanto, and Marc Hoyois for comments on a draft of this paper, and an anonymous referee for suggesting that there should exist a simpler proof of our presentation lemma. This project has received funding from the NSERC Discovery grant RGPIN-2025-07114, ``Motivic cohomology: theory and applications'', from the European Research Council (ERC) under the European Union's Horizon 2020 research and innovation programme (grant agreement No. 101001474), and from the SFB 1085 Higher Invariants.

\section{Nisnevich-local presentation lemma over a base}\label{sectionpresentationlemma}

\vspace{-\parindent}
\hspace{\parindent}

Gabber's presentation lemma, which can be considered as an analogue in algebraic geometry of the tubular neighbourhood theorem from differential geometry, describes the geometry of smooth schemes over a field. In this section, our goal is to prove a Nisnevich-local variant of this lemma for smooth schemes over more general base schemes (Corollary~\ref{cor:presentation_lemma}). This result is the main geometric ingredient needed in Section~\ref{sectionGersten} to prove Gersten's injectivity (Theorem~\ref{theoremGersteninjectivity}). 

Gabber's presentation lemma \cite{gabber_gersten_conjecture} was originally stated over an infinite field. Using ideas from Poonen’s proof of Bertini’s theorem over finite fields, Hogadi and Kulkarni \cite{hogadi_kulkarni_gabber_finite_field} proved the same statement over a finite field. Schmidt and Strunk \cite{schmidt_strunk} then extended this result to a Nisnevich-local presentation lemma over discrete valuation rings with infinite residue fields. Their proof crucially relies on the results of Kai \cite{kai} on good compactifications over a discrete valuation ring with infinite residue field. The approach of Schmidt and Strunk was later extended to a presentation lemma over noetherian domains with infinite residue fields by Deshmukh--Hogadi--Kulkarni--Yadav \cite{deshmukh-hogadi-kulkarni-yadav_presentation_lemma_noetherian_domains}. Independently, a weaker form of their Nisnevich-local presentation lemma, where the important finiteness condition is not guaranteed, was proved by Druzhinin \cite{druzhinin_a1_connectivity} over general noetherian rings.

Our proof of the presentation lemma is independent from the approaches taken in \cite{schmidt_strunk,deshmukh-hogadi-kulkarni-yadav_presentation_lemma_noetherian_domains,druzhinin_a1_connectivity}.\footnote{In the first version of this article, we gave a different proof of this presentation lemma over a general base scheme, which did not rely on the known presentation lemma over a field, and where we instead used the geometric techniques developed in \cite{ces_grothendieck-serre,gabber-liu-lorenzini} to extend Schmidt--Strunk's approach in terms of good compactifications to general base schemes.} We instead rely on the presentation lemma of Gabber and Hogadi--Kulkarni over a field, which we extend to general base schemes by a spreading out argument.

\begin{theorem}[Gabber \cite{gabber_gersten_conjecture}, Hogadi--Kulkarni \cite{hogadi_kulkarni_gabber_finite_field}]\label{theorempresentationlemmaoverafield}
    Let $k$ be a field, $X$ be a smooth $k$-scheme of finite type, and $Z \hookrightarrow X$ be a closed immersion of $k$-schemes of positive codimension. Then for every point $x \in X$, there exist 
    \begin{enumerate}[label=(\alph*)]
        \item a Zariski neighbourhood $x \in X' \subseteq X$,
        \item an affine open subscheme $T \subseteq \mathbb{A}^{d-1}_{k}$ where $d := \emph{dim}(X)$,
        \item and an étale morphism $\psi\colon X'\to\mathbb{A}^1_T$
    \end{enumerate}
    such that, for $Z' \colonequals Z\times_X X'$, the morphism $\psi|_{Z'}\colon Z' \rightarrow \mathbb{A}^1_T$ is a closed immersion, the image $\psi(Z')$ is finite over $T$, and $Z' \cong \psi(Z') \times_{\mathbb{A}^1_T} X'$.
\end{theorem}

\begin{proof}
    If $k$ is an infinite field, this is \cite[Lemma~3.1]{gabber_gersten_conjecture}. If $k$ is a finite field, this is \cite[Theorem~1.1]{hogadi_kulkarni_gabber_finite_field}.
\end{proof}

\begin{corollary}[Presentation lemma]\label{cor:presentation_lemma}
    Let $S$ be a scheme, $X$ be a smooth finitely presented $S$-scheme, and $Z \hookrightarrow X$ be a closed immersion of $S$-schemes. Then for every point $x \in X$ lying over a point $s \in S$ such that $\codim_{X_s}(Z_s)>0$, there exist
    \begin{enumerate}[label=(\alph*)]
        \item affine Nisnevich neighbourhoods $s \in S' \rightarrow S$ and $x \in X' \rightarrow X \times_S S'$,
        \item a smooth $S'$-scheme $T$, 
        \item and an étale morphism $\psi\colon X'\to\mathbb{A}^1_T$
    \end{enumerate}
    such that, for $Z'\colonequals Z\times_X X'$, the morphism $\psi|_{Z'}\colon Z' \rightarrow \mathbb{A}^1_T$ is a closed immersion, the image $\psi(Z')$ is finite over $T$, and $Z' \cong \psi(Z') \times_{\mathbb{A}^1_T} X'$. In particular, the commutative diagram of $T$-schemes
    $$\begin{tikzcd}
        X' \setminus Z' \arrow[d]\arrow[r, hook] & X' \arrow[d, "\psi"]\\ 
        \mathbb{A}^1_T \setminus \psi(Z') \arrow[r, hook] & \mathbb{A}^1_T
    \end{tikzcd}$$
    is a Nisnevich square, {\it i.e.}, it is a pullback square of schemes, $\mathbb{A}^1_T \setminus \psi(Z') \hookrightarrow \mathbb{A}^1_T$ is an open immersion, $\psi$~is an étale morphism, and $\psi(Z') \times_{\mathbb{A}^1_T} X' \rightarrow \psi(Z')$ is an isomorphism.
\end{corollary}

\begin{proof}
    The special fibre $X_s$ is a noetherian scheme, so the vanishing ideal of $Z_s$ inside $X_s$ is finitely generated. By spreading out a finite set of generators of this vanishing ideal to a Zariski neighbourhood of $s \in S$, we can assume, up to replacing $Z$ by a bigger closed subscheme of $X$, that $Z$ is finitely presented over $S$. The claim is Nisnevich-local around $s \in S$, and all the involved information is of finite presentation, so we can assume that $S$ is a henselian local scheme with closed point $s$. The claim is moreover Zariski-local around $x \in X$, so we can assume that $X$ is affine.

    By Theorem~\ref{theorempresentationlemmaoverafield} applied to $Z_s \hookrightarrow X_s$, there exist a Zariski neighbourhood $x \in X'_s \subseteq X_s$, an affine open subscheme $T_s \subseteq \mathbb{A}^{d-1}_s$, and an étale morphism $\psi_s : X'_s \rightarrow \mathbb{A}^1_{T_s}$ such that, for $Z'_s := Z_s \times_{X_s} X'_s$, the morphism $\psi_s : Z'_s \rightarrow \mathbb{A}^1_{T_s}$ is a closed immersion, the image $\psi_s(Z'_s)$ is finite over $T_s$, and $Z'_s \cong \psi_s(Z'_s) \times_{\mathbb{A}^1_{T_s}} X'_s$.

    We now lift this presentation over the special fibre to the desired presentation over $S$. Let $t \in T_s$ be the image of $x \in X'_s$ under the composite map $X'_s \xrightarrow{\psi_s} \mathbb{A}^1_{T_s} \rightarrow T_s$. Up to replacing $T_s$ by a smaller Zariski neighbourhood of $t \in \mathbb{A}^{d-1}_s$, we can assume that $T_s=D(\overline{f})$ is a standard Zariski open of $\mathbb{A}^{d-1}_s$. Choosing a section $f$ of $\mathbb{A}^{d-1}_S$ lifting $\overline{f}$, let $T := D(f) \subseteq \mathbb{A}^{d-1}_S$ be the corresponding lift of $T_s$. Similarly, up to replacing $X'_s$ by a smaller Zariski neighbourhood of $x \in X_s$, there exists a Zariski neighbourhood $x \in X' \subseteq X$ whose special fibre is $X'_s$.

    To lift the morphism $\psi_s : X'_s \rightarrow \mathbb{A}^1_{T_s}$ to a morphism $\psi : X' \rightarrow \mathbb{A}^1_T$ of $S$-schemes, it suffices to produce morphisms of $S$-schemes $X' \rightarrow \mathbb{A}^1_S$ and $X' \rightarrow T$ reducing to the given morphisms on the special fibre. The scheme $X'$ is affine, so the map $\mathcal{O}(X') \rightarrow \mathcal{O}(X'_s)$ is surjective. We let $X' \rightarrow \mathbb{A}^1_S$ be the morphism obtained by lifting the given global section of $X'_s$ to a global section of $X'$. For the second morphism, there is similarly a lift $X' \rightarrow \mathbb{A}^{d-1}_S$ of the composite morphism $X'_s \rightarrow T_s \subseteq \mathbb{A}^{d-1}_s$. Up to replacing $X'$ by the smaller Zariski neighbourhood $X' \times_{\mathbb{A}^{d-1}_S} T$ of $x \in X$, this morphism induces the desired lift $X' \rightarrow T$ of $X_s \rightarrow T_s$. 

    The obtained morphism $\psi : X' \rightarrow \mathbb{A}^1_T$ is étale on a Zariski neighbourhood of $x \in X'$, as a consequence of the étaleness of $\psi_s : X'_s \rightarrow \mathbb{A}^1_{T_s}$. Indeed, $\psi$ is finitely presented by \stacks{02FV}, flatness is ensured by \stacks{039C}, $G$-unramifiedness is ensured by \stacks{02GF}, and étaleness then follows by \stacks{02GV}.

    Let $Z' \colonequals Z \times_X X'$. We shall first replace $X'$ by a Nisnevich neighbourhood of $x \in X'$ to ensure that $Z'$ is finite over $T$, and then make $\psi|_{Z'} : Z' \rightarrow \mathbb{A}^1_T$ into a closed immersion satisfying $Z' \cong \psi(Z') \times_{\mathbb{A}^1_T} X'$. On the special fibre, $Z'_s \cong \psi_s(Z'_s)$ is finite over $T_s$, and the locus where $Z' \rightarrow T$ is quasi-finite is open \stacks{01TI}. Replacing $X'$ by an affine Zariski neighbourhood of $x \in X'$ whose base change to $Z'$ is contained in this open locus, we can then assume that $Z' \rightarrow T$ is quasi-finite. By the structure theorem for quasi-finite algebras over a henselian local ring \cite[\href{https://stacks.math.columbia.edu/tag/04GG}{Tag 04GG}\,(13)]{stacks-project}, $Z' \times_T \text{Spec}(\mathcal{O}_{T,t}^h)$ is a disjoint union of a finite $\mathcal{O}_{T,t}^h$-scheme and of a scheme that is disjoint from the special fibre of $\mathcal{O}_{T,t}^h$. Since $x \in Z'_s$ lies in the finite locus and everything is of finite presentation, this decomposition spreads out to a Nisnevich neighbourhood of $t \in T$. Replacing $T$ by such a Nisnevich neighbourhood, and $X'$ by an open of the corresponding base change that avoids this second locus, we can then assume that $Z' \rightarrow T$ is finite. 

    We now prove that $\psi|_{Z'} : Z' \rightarrow \mathbb{A}^1_T$ is a closed immersion on a Zariski neighbourhood of $x \in X'$. First note that $\psi|_{Z'}$ is finite, since $\mathbb{A}^1_T \rightarrow T$ is separated and the composite $Z' \xrightarrow{\psi|_{Z'}} \mathbb{A}^1_T \rightarrow T$ is finite \stacks{035D}. Since $\psi|_{Z'}$ is a morphism of affine schemes, it suffices to prove that the map of finite type $\mathcal{O}_{\mathbb{A}^1_T}$-modules $f \colon \mathcal{O}_{\mathbb{A}^1_T} \rightarrow (\psi|_{Z'})_* \mathcal{O}_{Z'}$ is surjective on a Zariski neighbourhood of $x \in X'$. The base change of this map to the special fibre of $S$ is $\mathcal{O}_{\mathbb{A}^1_{T_s}} \rightarrow (\psi_s|_{Z'_s})_* \mathcal{O}_{Z'_s}$ \stacks{01LX}, which is surjective because $\psi_s|_{Z'_s}$ is a closed immersion. In particular, the cokernel of $f$ is a finite $T$-module whose support is disjoint from $T_s$. By Nakayama's lemma, this implies that this cokernel vanishes on a Zariski neighbourhood of $t \in T$. Base changing $X'$ and $Z'$ to this Zariski neighbourhood of $t \in T$, we can assume that $\psi|_{Z'}$ is a closed immersion. In particular, we get that $\psi|_{Z'} : Z' \xrightarrow{\sim} \psi(Z')$ and that $\psi(Z')$ is finite over $T$.

    Finally, we ensure that $Z' \cong \psi(Z') \times_{\mathbb{A}^1_T} X'$. By construction, we already know that $Z' \subseteq \psi^{-1}(\psi(Z'))$, and that the base change $\psi^{-1}(\psi(Z')) = \psi(Z') \times_{\mathbb{A}^1_T} X' \rightarrow \psi(Z')$ of $\psi$ is \'etale. The composite map $Z' \rightarrow \psi(Z')$ is moreover an isomorphism by the previous paragraph, so $Z'$ is open and closed in $\psi^{-1}(\psi(Z'))$. Removing the complementary component (which is closed and does not contain $x \in Z'$) from $X'$, we obtain $Z' = \psi(Z') \times_{\mathbb{A}^1_T} X'$, which concludes the proof.
\end{proof}

\begin{remark}
    The main difference between our result and Gabber's presentation lemma over a field is that our presentation lemma is Nisnevich-local rather than Zariski-local. All the sheaves considered in this article being in fact Nisnevich sheaves, this difference will not be important for our purposes. Upgrading the results of this section to Zariski-local statements could be interesting for other applications, but it would seem to require new ideas.
\end{remark}

\begin{remark}
    When $S$ is the spectrum of a Dedekind domain with infinite residue fields, the previous result was proved by Schmidt--Strunk (\ccite{schmidt_strunk}{Theorem~2.4}). To the best of our knowledge, the previous result is new when $S$ is the spectrum of a Dedekind domain with not necessarily infinite residue fields ({\it e.g.}, when $S=\text{Spec}(\Z)$).
\end{remark}

\section{Gersten's injectivity in mixed characteristic}
\label{sectionGersten}

\vspace{-\parindent}
\hspace{\parindent}

In this section, we apply the presentation lemma of the previous section to prove a Nisnevich-local version of Gersten's injectivity for smooth schemes over not necessarily discrete valuation rings (Corollary~\ref{cor:theoremGersteninjectivity}). This result applies in particular for $\mathbb{A}^1$-invariant Nisnevich sheaves (Remark~\ref{remarkA1invimpliesdeflatable}) and for non-$\mathbb{A}^1$-invariant motivic spectra in the sense of \cite{annala_atiyah_2024} (Corollary~\ref{cor:Gersten_injectity_motivic_spectrum_annala_iwasa}), and will be the main new ingredient in the proofs of Theorems~\ref{theoremintrocomparisonclassicalmotivic} and \ref{theoremintroBL}.

We first prove a Gersten's injectivity statement in the context of Nisnevich sheaves of pointed sets
(Theorem~\ref{theoremGersteninjectivity}). The main property that we ask for these Nisnevich sheaves is the Horrocks principle (Definition~\ref{defn:generalised_horrocks}), whose terminology is inspired by the literature on $G$-torsors (see for instance \cite[Proposition~$2.1.5$]{ces_torsors}), and which is a non-commutative variant of the axiom \say{SUB2} of \cite{colliot-thelene_bloch_1997} (see Definition~\ref{definitiondeflatablepresheaves}). We then apply this result to the homotopy groups of spectra-valued Nisnevich sheaves to prove our main Gersten's injectivity result (Corollary~\ref{cor:theoremGersteninjectivity}).

For the rest of this section, given a scheme $X$, we denote by $\infty_X$, $j_X$, $\pi_X$ and $\bar{\pi}_X$ the canonical morphisms of schemes appearing in the following commutative diagram.
\begin{equation}\label{diag:notations_gersten_injectivity}\tag{$\ast$}
    \begin{tikzcd}\mathbb{A}^1_{X}\arrow{dr}[swap]{\pi_{X}}\arrow[r,hook,"j_{X}"]&\mathbb{P}^1_{X}\arrow[d, "\bar{\pi}_X"]&X\arrow[hook']{l}[swap]{\infty_{X}}\arrow[ld,-,double equal sign distance,double]\\ &X\end{tikzcd}
\end{equation}

\begin{definition}[Horrocks principle]\label{defn:generalised_horrocks}
    Given a scheme $S$, a functor $F\colon \text{Sch}^{\text{qcqs,op}}_S \rightarrow \text{Set}_{\ast}$ 
    is said to satisfy the \textit{Horrocks principle} if for any henselian local $S$-scheme $X$, every class $[\sigma]\in F(\mathbb{P}^1_X)$ that trivialises at $\infty$ ({\it i.e.}, $\infty_X^{\ast}[\sigma]=[\ast]$) also trivialises on $\mathbb{A}^1_X$ ({\it i.e.},~$j_X^{\ast}[\sigma]=[\ast])$.    
\end{definition}

Our argument to establish Gersten's injectivity follows the standard techniques axiomatised by Colliot-Thélène--Hoobler--Kahn \cite{colliot-thelene_bloch_1997}, where we replace the reduction to the case of fields by a reduction to the case of henselian valuation rings. The following result is a crucial step in proving our main Gersten's injectivity result (Theorem~\ref{theoremGersteninjectivity}), and is the only point in our argument where we use the presentation lemma (Corollary~\ref{cor:presentation_lemma}).

\begin{theorem}\label{Nisnevich_local_Gersten_Injectivity}
    Let $P$ be a Prüfer ring,\footnote{Recall that a commutative ring is a \textit{Prüfer domain} if it is an integral domain whose localisation at every prime ideal is a valuation ring, and is a \textit{Prüfer ring} if it is a finite product of Prüfer domains. In particular, a commutative ring is a valuation ring if and only if it is a local Prüfer ring.} $R$ be the henselisation of an essentially smooth local $P$-algebra, and $s \in \spec(P)$ be the image of the closed point of $\spec(R)$. Then the fibre $\spec(R)_s$ of $\spec(R)$ over $s$ is irreducible. Moreover, for every finitary Nisnevich sheaf $F\colon \emph{Sch}_P^{\emph{qcqs,op}} \rightarrow \emph{Set}_\ast$ satisfying the Horrocks principle in the sense of Definition~\ref{defn:generalised_horrocks}, we have that
    $$\ker(F(R)\to F(R_{\eta}))=\{\ast\},$$
    where $\eta\in\spec(R)$ is the generic point of the irreducible scheme $\spec(R)_s$.
\end{theorem}

\begin{proof}
    The claims of interest are insensitive to replacing the Prüfer ring $P$ by the henselisation of its localisation at $s$. In particular, we can assume that $P$ is a henselian 
    valuation ring with maximal ideal $s$. The scheme $\spec(R)$ is local, thus, so is the scheme $\spec(R)_s$. The essential smoothness hypothesis on the $P$-algebra $R$ ensures that the scheme $\text{Spec}(R)_s$ is moreover regular. Consequently, the scheme $\spec(R)_s$ is irreducible (\stacks{00NP}). To prove the second claim, because the presheaf $F$ commutes with filtered colimits of rings, we can further reduce to the case where the valuation ring $P$ is of finite rank (see for instance \ccite{amar:unramified_case_gersten_conjecture}{Lemma~$2.5\,(b)$}).

    Let $[\sigma] \in F(R)$ be a class whose image in $F(R_\eta)$ vanishes. We want to prove that $[\sigma]$ itself vanishes.
    For the rest of this proof, given a closed immersion of $P$-schemes $Z \hookrightarrow X$, we denote by $F_Z(X)$ the pointed set $\text{ker}(F(X)\to F(X \setminus Z))$. The presheaf $F$ commutes with filtered colimits of rings, so there exist a smooth affine finitely presented $P$-scheme $X$ and a closed immersion of $P$-schemes $Z \hookrightarrow X$ such that $R$ is the henselisation of the local ring of $X$ at some point $x \in X$ lying over $s \in \text{Spec}(P)$ and such that the class $[\sigma] \in F(R)$ is the image of a class $[\sigma] \in F(X)$ that vanishes in $F(X \setminus Z)$, {\it i.e.}, $[\sigma] \in F(R)$ is the image of a class $[\sigma] \in F_Z(X)$.

    To prove the desired result, it then suffices to produce a Nisnevich neighbourhood $x \in X' \rightarrow X$ which satisfies, for $Z' \colonequals Z \times_X X' \hookrightarrow X'$, that the pullback class $[\sigma'] \in F_{Z'}(X')$ vanishes. If $R$ is of relative dimension zero over $P$, then the local ring $\mathcal{O}_{X,x}$ is étale over $P$, and is even a valuation ring with maximal ideal $x$ (\stacks{0ASJ}). In particular, the local ring $R$ has maximal ideal $\eta$, and $R=R_{\eta}$. We assume now that $R$ is of positive relative dimension over $P$. By \ccite{amar:unramified_case_gersten_conjecture}{Lemma 2.13\,(2)}, where we use that the valuation ring $P$ is of finite rank, we can then further assume that $\dim(Z_s) < \dim(X_s)$.
    
    Therefore, by the presentation lemma proved in the previous section over the henselian local scheme~$S$ (Corollary~\ref{cor:presentation_lemma}), there exist a smooth $S$-scheme~$T$, a Nisnevich neighbourhood $X' \rightarrow X$ of $x$, and a Nisnevich square
    $$\begin{tikzcd}
        X' \setminus Z' \arrow[r,hook] \ar[d] & X' \arrow[d,"\psi"] \\
        \mathbb{A}^1_T \setminus \psi(Z') \arrow[r,hook] & \mathbb{A}^1_T
    \end{tikzcd}$$
    such that the induced morphism $Z' \xrightarrow{\cong} \psi(Z') \rightarrow T$ is finite. Let $t \in T$ be the image of $x \in X'$ via the composite morphism $X' \xrightarrow{\psi} \mathbb{A}^1_T \rightarrow T$. By Nisnevich excision, the natural map $F_{\psi(Z')}(\mathbb{A}^1_T) \rightarrow F_{Z'}(X')$ is an isomorphism, and it then suffices to prove that the class $[\sigma'] \in F_{\psi(Z')}(\mathbb{A}^1_T)$ vanishes Nisnevich-locally around $t \in T$. By finiteness of the map $\psi(Z') \rightarrow T$, the composite morphism $\psi(Z') \hookrightarrow \mathbb{A}^1_T \hookrightarrow \mathbb{P}^1_T$ is a closed immersion (\cite[\href{https://stacks.math.columbia.edu/tag/01WN}{Tags 01WN} and \href{https://stacks.math.columbia.edu/tag/01W6}{01W6}]{stacks-project}), hence the commutative diagram of $T$-schemes
    $$\begin{tikzcd}
        \mathbb{A}^1_T \setminus \psi(Z') \arrow[r,hook] \ar[d] & \mathbb{A}^1_T \arrow[d,"j_T"] \\
        \mathbb{P}^1_T \setminus \psi(Z') \arrow[r,hook] & \mathbb{P}^1_T
    \end{tikzcd}$$
    is a Nisnevich square. Again, by Nisnevich excision, the natural map $F_{\psi(Z')}(\mathbb{P}^1_T) \rightarrow F_{\psi(Z')}(\mathbb{A}^1_T)$ is an isomorphism, hence there exists a class $[\tilde{\sigma}'] \in F(\mathbb{P}^1_T)$ whose image in $F(\mathbb{P}^1_T \setminus \psi(Z'))$ vanishes and such that $j_T^\ast [\tilde{\sigma}'] = [\sigma'] \in F(\mathbb{A}^1_T)$. By construction, the closed subscheme $\psi(Z') \subseteq \mathbb{P}^1_T$ does not meet the $\infty$\nobreakdash-section of $\mathbb{P}^1_T$, so the map $\infty_T^\ast\colon F(\mathbb{P}^1_T) \rightarrow F(T)$ naturally factors through the map $F(\mathbb{P}^1_T) \rightarrow F(\mathbb{P}^1_T \setminus \psi(Z'))$. In particular the class $\infty_T^\ast [\tilde{\sigma}'] \in F(T)$ vanishes. By the Horrocks principle (Definition~\ref{defn:generalised_horrocks}), the class $j_T^{\ast}[\tilde{\sigma}'] \in F(\mathbb{A}^1_T)$ then also vanishes Nisnevich-locally around $t \in T$, which concludes the proof.
\end{proof}

The final ingredient in the proof of Theorem~\ref{theoremGersteninjectivity} is Corollary~\ref{prop:vansishing_valuation_rings}, which is obtained as a consequence of the following result.
\begin{lemma}\label{lem:gersten_injectivity_valuation_ring}
    Let $\ca{V}$ be a subcategory of the category of commutative rings that is closed under localisations of rings at multiplicative subsets, and $F\colon \ca{V} \rightarrow \emph{Set}_{\ast}$ be a finitary Nisnevich sheaf. Then, the following are equivalent:
    \begin{enumerate}
        \item\label{lem:gersten:pt:1} for every henselian valuation ring $V \in \ca{V}$ of finite rank, we have
        $$\ker(F(V)\to F(\emph{Frac}(V)))=\{\ast\};$$
        \item\label{lem:gersten:pt:2} for every semilocal Prüfer domain $P \in \ca{V}$, we have
        $$\ker(F(P)\to F(\emph{Frac}(P)))=\{\ast\}.$$
    \end{enumerate}
\end{lemma}
\begin{proof}
    Every henselian valuation ring is a semilocal Prüfer domain, so it suffices to prove that \eqref{lem:gersten:pt:1} $\Rightarrow$ \eqref{lem:gersten:pt:2}. Let $P$ be a semilocal Prüfer domain. The presheaf $F$ being finitary, we may assume by a limit argument that the Prüfer domain $P$ is of finite Krull dimension (\ccite{amar:unramified_case_gersten_conjecture}{Lemma~$2.5$}). Letting $\text{MaxSpec}(P)$ be the maximal spectrum of $P$, we prove the desired result by induction on the integer
    \begin{equation*}
        d=\delta(P)\colonequals \sum_{\mathfrak{m}\in\text{MaxSpec}(P)}\text{dim}(P_{\mathfrak{m}})<\infty.
    \end{equation*} If $d=0$, then $P=\text{Frac}(P)$, and the claim is true. We assume now that $d>0$, and that the result is true for all semilocal Prüfer domains $P'$ satisfying $\delta(P')<d$. Let $\mathfrak{m}$ be a maximal ideal of $P$, and let $a$ be an element of $P$ whose vanishing locus equals the subset $\{\mathfrak{m}\}\subseteq\spec(P)$.\footnote{Since $P$ has finitely many prime ideals, such an $a$ always exist by prime avoidance \stacks{00DS}.} The ring $P[\frac{1}{a}]$ then satisfies $\delta(P[\tfrac{1}{a}])=d-1$ (since $\text{Spec}(P[\frac{1}{a}])=\text{Spec}(P)\setminus\{\mathfrak{m}\}$), and is a semilocal Prüfer domain (\ccite{amar:unramified_case_gersten_conjecture}{Lemma~$2.2\,(1)$}). Consequently, by the induction hypothesis, it suffices to prove the equality
    $$\ker(F(P)\to F(P[\tfrac{1}{a}]))=\{\ast\}.$$
    By Nisnevich excision, the left-hand side of this equality is naturally identified with the pointed set $\ker(F(P^h_{\mathfrak{m}})\to F(P^h_{\mathfrak{m}}[\frac{1}{a}]))$, where $P^h_{\mathfrak{m}}$ is the henselisation of $P$ along the ideal $\mathfrak{m}$. The henselisation of a local ring of a Prüfer domain is a henselian valuation ring, so the result is a consequence of $(1)$.
\end{proof}

\begin{corollary}\label{prop:vansishing_valuation_rings}
    Let $\cal{V}$ be a subcategory of the category of commutative rings that is closed under localisations of rings at multiplicative subsets, and $F\colon\cal{V}\to\emph{Set}_{\ast}$ be a finitary Nisnevich sheaf. If $F(V)=\{\ast\}$ for every henselian valuation ring $V \in \cal{V}$ of finite rank, then $F(P)=\{\ast\}$ for every semilocal Prüfer domain $P \in \cal{V}$.
\end{corollary}

\begin{proof}
    Let $P \in \cal{V}$ be a semilocal Prüfer domain. Any field is a henselian valuation ring of finite rank, so that $F(\text{Frac}(P)) = \{\ast\}$, and in particular $F(P)=\ker(F(P)\rightarrow F(\Frac(P))$. Lemma~\ref{lem:gersten_injectivity_valuation_ring} and the assumption on henselian valuation rings of finite rank then imply the desired result.
\end{proof}

\begin{theorem}\label{theoremGersteninjectivity}
    Let $P$ be a Prüfer ring and $F\colon \emph{Sch}^{\emph{qcqs,op}}_P\rightarrow\emph{Set}_{\ast}$ be a finitary Nisnevich sheaf satisfying the following conditions:
    \begin{enumerate}[label=(\roman*)]
        \item the presheaf $F$ satisfies the Horrocks principle in the sense of Definition \ref{defn:generalised_horrocks};
        \item\label{pt:2:gersten_injectivity} for every henselian valuation ring $V$ which is an ind-smooth $P$-algebra, we have $F(V)=\{\ast\}.$
    \end{enumerate}
    Then, for every ind-smooth $P$-scheme $X$, we have $F(X)=\{\ast\}$.
\end{theorem}
\begin{proof}
    The presheaf $F$ being a finitary Nisnevich sheaf, we can assume that $X=\text{Spec}(R)$, where $R$ is the henselisation of an essentially smooth local $P$-algebra, and that $P$ is a valuation ring of finite rank (see for instance \ccite{amar:unramified_case_gersten_conjecture}{Lemma~$2.5$}). By Theorem~\ref{Nisnevich_local_Gersten_Injectivity}, the $P$-fibre of the closed point of $\text{Spec}(R)$ admits a unique generic point $\eta \in \text{Spec}(R)$, and the localisation $A\colonequals R_{\eta}$ satisfies that
    $$\ker(F(R)\to F(A))=\{\ast\},$$
    hence it suffices to prove that $F(A)=\{\ast\}$. The local $P$-algebra $A$ is moreover a semilocal Prüfer domain of finite Krull dimension (\ccite{amar:unramified_case_gersten_conjecture}{Lemma~$2.13\,(1)$}).
    The hypothesis (ii) and Corollary~\ref{prop:vansishing_valuation_rings} then imply that $F(A)=\{\ast\}$.
\end{proof}

Although Theorem~\ref{theoremGersteninjectivity} is stated in the generality of presheaves of pointed sets, our interest in this article lies primarily in presheaves of spectra. Accordingly, in the rest of this section, we explain how this result can be applied to certain spectra-valued presheaves, which we call {\it deflatable} (Definition~\ref{definitiondeflatablepresheaves}). This condition essentially corresponds to the axiom \say{SUB2} of \cite{colliot-thelene_bloch_1997}, and is equivalent to the one introduced in \cite[Definition~$6.7$]{elmanto_motivic_2023}. Note that any $\mathbb{A}^1$-invariant presheaf is deflatable (Remark~\ref{remarkA1invimpliesdeflatable}), and so is any family of functors satisfying the $\mathbb{P}^1$-bundle formula (Lemma~\ref{lemmapbfimpliesdeflatable}).

\begin{definition}[Deflatability]\label{definitiondeflatablepresheaves}
    Given a scheme $S$, a functor $\scr{F}\colon \text{Sch}^{\text{qcqs,op}}_S \rightarrow \text{Sp}$ is called \emph{deflatable} if for any qcqs $S$-scheme $X$, there exists a functorial equivalence between the maps
    $$\begin{tikzcd}
        \scr{F}(\mathbb{P}^1_X) \arrow[shift right=1.25ex]{rr}[swap]{(\infty_X\circ\pi_X)^{\ast}}\arrow[rr, shift left=1.25ex, "j_X^{\ast}"]&& \scr{F}(\mathbb{A}^1_X),
    \end{tikzcd}$$
    where the maps $\infty_X$, $j_X$, and $\pi_X$ are defined in \eqref{diag:notations_gersten_injectivity}.
\end{definition}

\begin{remark}
    In the previous definition, note that a functor $\scr{F}\colon\text{Sch}^{\text{qcqs,op}}_S \rightarrow \text{Sp}$ is deflatable as soon as the category of equivalences $\text{Eq}(j^{\ast},(\infty\circ\pi)^{\ast})$ between the functors $j^\ast$ and $(\infty \circ \pi)^{\ast}$ is nonempty. In particular, any such equivalence is not part of the data (see also \cite[Remark~$6.8$]{elmanto_motivic_2023}).
\end{remark}

\begin{remark}\label{remarkdeflatableimpliesHorrocks}
    Let $S$ be a scheme, and $\scr{F}\colon\text{Sch}^{\text{qcqs,op}}_S \rightarrow \text{Sp}$ be a deflatable functor. Then for every integer $i \in \Z$, the functor $\pi_i(\scr{F}(-))\colon \text{Sch}^{\text{qcqs,op}}_S \rightarrow \text{Set}_{\ast}$ satisfies the Horrocks principle in the sense of Definition~\ref{defn:generalised_horrocks}.
\end{remark}

\begin{remark}\label{remarkA1invimpliesdeflatable}
    Let $S$ be a qcqs scheme, and $\scr{F}\colon\text{Sch}^{\text{qcqs,op}}_S \rightarrow \text{Sp}$ be a functor. If $\scr{F}$ is $\mathbb{A}^1$-invariant, then $\scr{F}$ is deflatable. Indeed, in this case the natural map $\pi^\ast$ is an equivalence and, for any section $s \in \mathbb{P}^1_S(S)\setminus\{\infty_S\}$, both $j^\ast$ and $\infty^\ast$ are $\mathbb{A}^1$-homotopy equivalent to $s^\ast$.
\end{remark}

\begin{definition}[$\mathbb{P}^1$-bundle formula]\label{definitionP^1bundleformula}
    Let $S$ be a scheme. A family of presheaves $\scr{F}(i) : \text{Sch}^{\text{qcqs,op}}_S \rightarrow \text{Sp}$, $i \in \Z$, equipped for every $i \in \Z$ with a map of functors $c_1 : \text{Pic}(-) \rightarrow \text{Hom}(\scr{F}(i-1)(-)[-2],\scr{F}(i)(-))$, is said to {\it satisfy the $\mathbb{P}^1$-bundle formula} if for every qcqs $S$-scheme $X$ and every integer $i \in \Z$, the natural map
    $$\bar{\pi}^\ast_X \oplus c_1(\mathcal{O}(1))\, \bar{\pi}_X^\ast : \scr{F}(i)(X) \oplus \scr{F}(i-1)(X)[-2] \longrightarrow \scr{F}(i)(\mathbb{P}^1_X)$$
    is an equivalence of spectra.
\end{definition}

\begin{remark}
    In all the situations of interest, the map $c_1$ in Definition~\ref{definitionP^1bundleformula} is induced by a suitable notion of first Chern class, corresponding to an orientation for the family of presheaves $\{\scr{F}(i)\}_{i \in \Z}$ (\cite[Section~$3$]{annala_motivic_2023}).
\end{remark}

\begin{lemma}[\cite{colliot-thelene_bloch_1997}]\label{lemmapbfimpliesdeflatable}
Let $S$ be a qcqs scheme, and $\scr{F}(i)$, $i \geqslant 0$, be spectra-valued presheaves on smooth $S$\nobreakdash-schemes (resp. on qcqs $S$-schemes). If the family of presheaves $\{\scr{F}(i)\}_{i \geqslant 0}$ satisfies the $\mathbb{P}^1$-bundle formula, then for every integer $i \geqslant 0$, the presheaf $\scr{F}(i)$ is deflatable.
\end{lemma}

\begin{proof}
    The proof is the same as in \cite[Proposition~$5.4.3$]{colliot-thelene_bloch_1997} (see also \cite[Lemma~$6.12$]{elmanto_motivic_2023}). More precisely, consider the diagram of spectra
    \begin{equation*}
        \begin{tikzcd}
            \scr{F}(i)(\bb{P}^1_X)\arrow[rr,"j_X^{\ast}"]\arrow[drr,"\infty_X^{\ast}"]&&F(i)(\bb{A}^1_X)\\\scr{F}(i)(X)\oplus \scr{F}(i-1)(X)[-2]\arrow[u, "\sim"{swap, sloped}]
            &&\scr{F}(i)(X),\arrow[u,"\pi_X^{\ast}"]
        \end{tikzcd}
    \end{equation*}
    where the left map is an equivalence by the $\mathbb{P}^1$-bundle formula. We want to prove that the right triangle commutes. It suffices to do so after pre-composing with the left equivalence, and to check the commutativity on each of the factors. On the first factor, the commutativity follows by functoriality of $\scr{F}(i)$. On the second factor, it is a consequence of the functoriality of the map $c_1$, and of the identifications $j^{\ast}_X(\ca{O}_{\bb{P}^1_X}(1))\cong \ca{O}_{\bb{A}^1_X}$ and $\infty_X^{\ast}(\ca{O}_{\bb{P}^1_X}(1))\cong \ca{O}_X$.
\end{proof}

\begin{corollary}\label{cor:theoremGersteninjectivity}
    Let $P$ be a Prüfer ring ({\it e.g.}, a valuation ring), $\scr{F}\colon \emph{Sch}^{\emph{qcqs,op}}_P \rightarrow \emph{Sp}$ be a finitary Nisnevich sheaf, and $j \in \Z$ be an integer satisfying the following conditions:
    \begin{enumerate}[label=(\roman*)]
        \item the presheaf $\scr{F}$ is deflatable in the sense of Definition~\ref{definitiondeflatablepresheaves};
        \item for any henselian valuation ring $V$ which is an ind-smooth $P$-algebra, we have $\pi_j(\scr{F}(V))=0.$
    \end{enumerate}
    Then, for any henselian local ind-smooth $P$-algebra $R$, we have $\pi_j(\scr{F}(R))=0$.
\end{corollary}

\begin{proof}
    This is a consequence of Theorem~\ref{theoremGersteninjectivity} where we take $F$ to be the Nisnevich sheafification of the presheaf $\pi_j(\mathscr{F}(-))$, which satisfies the Horrocks principle by Remark~\ref{remarkdeflatableimpliesHorrocks}.
\end{proof}

\begin{remark}\label{rem:gersten_injectivity_noetherian}
    If $P$ is a field (resp. a discrete valuation ring), then the local ring $R_{\eta}$ appearing in Theorem~\ref{Nisnevich_local_Gersten_Injectivity} is also a field (resp. a discrete valuation ring). As a consequence, the results of this section actually extend the results of Colliot-Thélène--Hoobler--Khan \cite{colliot-thelene_bloch_1997} (resp. of Gillet--Levine \cite{gillet_levine}).
\end{remark}
\begin{corollary}\label{cor:Gersten_injectity_motivic_spectrum_annala_iwasa}
    Let $P$ be a Prüfer ring, and $\{E(i)\}_{i \in \Z} \in \emph{MS}(P)$ be a $\mathbb{P}^1$-motivic spectrum in the sense of \cite{annala_atiyah_2024}. Then for all integers $i,j \in \Z$, if the cohomology presheaf $\emph{H}^j(E(i)(-))$ vanishes on henselian valuation rings which are ind-smooth $P$-algebras, then it vanishes on all henselian local ind-smooth $P$\nobreakdash-algebras.
\end{corollary}

\begin{proof}
    The presheaf $\text{H}^j(E(i)(-))$ is deflatable by \cite[Corollaries~$4.4$ and~$4.11$]{annala_algebraic_2025} (for oriented $\mathbb{P}^1$-motivic spectra, this is more directly a consequence of Lemma~\ref{lemmapbfimpliesdeflatable}), so the result is a consequence of Corollary~\ref{cor:theoremGersteninjectivity}.
\end{proof}

\section{Review of non-$\bb{A}^1$-invariant motivic cohomology}
\label{sectionreview}

\vspace{-\parindent}
\hspace{\parindent}

Given a qcqs scheme $X$, one may consider the associated spectra given by non-connective algebraic $K$-theory $\text{K}(X)$, topological cyclic homology $\text{TC}(X)$, Weibel's homotopy $K$-theory $\text{KH}(X) \colonequals (L_{\mathbb{A}^1} \text{K})(X)$, and cdh-local topological cyclic homology $(L_{\text{cdh}} \text{TC})(X)$. By a theorem of Kerz--Strunk--Tamme \cite{kerz_algebraic_2018} and Land--Tamme \cite{land_k-theory_2019}, these invariants are related by a functorial cartesian square
$$\begin{tikzcd}
    \text{K}(X)\arrow[d]\arrow[r]&\text{TC}(X)\arrow[d]\\ \text{KH}(X)\arrow[r]&(L_{\text{cdh}}\text{TC})(X),
\end{tikzcd}$$
where the top horizontal map is the cyclotomic trace map, and where the vertical maps are the canonical maps. One may then consider the following {\it motivic filtrations} on these invariants.

\begin{enumerate}
    \item The $\N$-indexed filtration $\text{Fil}^\star_{\text{cdh}} \text{KH}(X)$, constructed by Bachmann--Elmanto--Morrow \cite{bachmann_A^1-invariant_2025} in terms of the classical $\mathbb{A}^1$-invariant motivic filtration on the algebraic $K$-theory of smooth affine $\Z$-schemes. We denote by $\Z(i)^{\text{cdh}}(X) \colonequals \text{gr}^i_{\text{cdh}} \text{KH}(X)[-2i]$ the shifted graded pieces of this filtration, which provide a theory of cdh-local motivic cohomology for $X$, where the cdh-topology is defined as in \cite{ehik}.
    \item The $\Z$-indexed filtration $\text{Fil}^\star_{\text{mot}} \text{TC}(X)$, constructed by Elmanto--Morrow \cite{elmanto_motivic_2023} and the first-named author \cite{bouis_motivic_2024} in terms of the Hochschild--Kostant--Rosenberg filtration on negative cyclic homology \cite{antieau_periodic_2019} and of the Bhatt--Morrow--Scholze filtration on $p$-completed topological cyclic homology \cite{bhatt_topological_2019} (or, alternatively, in terms of the even filtration \cite{hahn_motivic_2022}). We denote by $\Z(i)^{\text{TC}}(X) \colonequals \text{gr}^i_{\text{mot}}\text{TC}(X)[-2i]$ the shifted graded pieces of this filtration, which recover the Hodge-completed derived de Rham cohomology of $X$ in characteristic zero, and the syntomic cohomology of $X$ in positive characteristic.
    \item The $\Z$-indexed filtration $\text{Fil}^\star_{\text{mot}} L_{\text{cdh}} \text{TC}(X)$, introduced by Elmanto--Morrow \cite{elmanto_motivic_2023} and the first-named author \cite{bouis_motivic_2024} as the cdh sheafification $(L_{\text{cdh}} \text{Fil}^\star_{\text{mot}} \text{TC})(X)$ of the previous filtration. Hense\-lian valuation rings being the local rings for the cdh topology, the structural properties of the shifted graded pieces $L_{\text{cdh}} \Z(i)^{\text{TC}}(X)$ are usually determined by the value of $\Z(i)^{\text{TC}}$ at henselian valuation rings.
\end{enumerate}

These three filtrations are multiplicative ({\it i.e.}, have a structure of filtered $\mathbb{E}_{\infty}$-rings), and are functorial in~$X$. Given these filtrations, Elmanto--Morrow \cite{elmanto_motivic_2023} and the first-named author \cite{bouis_motivic_2024} then proved that there exists a multiplicative, functorial, $\N$-indexed filtration $\mathrm{Fil}^\star_{\text{mot}}\mathrm{K}(X)$ on $\text{K}(X)$ such that the previous cartesian square upgrades to a filtered cartesian square. We denote by $\Z(i)^{\text{mot}}(X) \colonequals \text{gr}^i_{\text{mot}}\text{K}(X)[-2i]$ the shifted graded pieces of this filtration, which provide a theory of motivic cohomology for $X$. By construction, the motivic complexes $\Z(i)^{\text{mot}}(X)$ are related to the previous invariants by a functorial cartesian square

$$\begin{tikzcd}
    \Z(i)^{\text{mot}}(X) \ar[r] \ar[d] & \Z(i)^{\text{TC}}(X) \ar[d] \\
    \Z(i)^{\text{cdh}}(X) \ar[r] & \big(L_{\text{cdh}} \Z(i)^{\text{TC}}\big)(X)
\end{tikzcd}$$
in the derived category $\mathcal{D}(\Z)$. In the rest of this section, we review the results about the motivic complexes $\Z(i)^{\text{mot}}(X)$ that will be used in later sections.

\begin{theorem}[\cite{elmanto_motivic_2023,bouis_motivic_2024,bouis_weibel_2025}]\label{theoremreview}
    The functors
    $$\Z(i)^{\emph{mot}}\colon\emph{Sch}^{\emph{qcqs,op}} \longrightarrow \mathcal{D}(\Z), \quad i \geqslant 0$$
    are finitary Nisnevich sheaves (\cite[Corollary~$4.59$]{bouis_motivic_2024}), and satisfy the following properties for every qcqs scheme $X$.
    \begin{enumerate}
        \item\label{pt:P^1_bundle_formula} The $\mathbb{P}^1$-bundle formula, {\it i.e.}, the natural maps
        $$\pi^\ast \oplus c_1^{\emph{mot}}(\mathcal{O}(1)) \pi^\ast\colon\Z(i)^{\emph{mot}}(X) \oplus \Z(i-1)^{\emph{mot}}(X)[-2] \longrightarrow \Z(i)^{\emph{mot}}(\mathbb{P}^1_X), \quad i \geqslant 0,$$
        where $\pi\colon\mathbb{P}^1_X \rightarrow X$ is the canonical projection map, are equivalences in the derived category $\mathcal{D}(\Z)$ (\cite[Theorem~$4.7$]{bouis_weibel_2025}). More generally, the functors $\Z(i)^{\emph{mot}}$ satisfy the projective bundle formula (\cite[Theorem~$4.19$]{bouis_weibel_2025}).
        \item There exists a natural equivalence
        $$\emph{K}(X)_{\Q} \simeq \bigoplus_{i \geqslant 0} \Q(i)^{\emph{mot}}(X)[2i]$$
        in the derived category $\mathcal{D}(\Q)$, induced by suitable Adams operations on the rational algebraic $K$-theory $\emph{K}(X)_{\Q}$ (\cite[Corollary~$4.60$]{bouis_motivic_2024}). 
        \item\label{pt:syntomic_period_map} For every prime number $p$ and every integer $k \geqslant 1$, there is a natural fibre sequence
        $$\Z/p^k(i)^{\emph{mot}}(X) \longrightarrow \Z/p^k(i)^{\emph{syn}}(X) \longrightarrow \big(L_{\emph{cdh}} \tau^{>i} \Z/p^k(i)^{\emph{syn}}\big)(X)$$
        in the derived category $\mathcal{D}(\Z/p^k)$ (\cite[Theorem~$5.10$]{bouis_motivic_2024}). In particular, the left map of this fibre sequence is an isomorphism in degrees less than or equal to $i$ \cite[Corollary~$5.11$]{bouis_motivic_2024}.
    \end{enumerate}
\end{theorem}

In Theorem~\ref{theoremreview}\,$(3)$, we use the {\it syntomic realisation map} 
$$\Z/p^k(i)^{\text{mot}}(X) \longrightarrow \Z/p^k(i)^{\text{syn}}(X)$$
from the motivic complex $\Z/p^k(i)^{\text{mot}}(X) \colonequals \Z(i)^{\text{mot}}(X) \otimes_{\Z}^{\mathbb{L}} \,\Z/p^k$ to Bhatt--Lurie's syntomic complex $\Z/p^k(i)^{\text{syn}}(X)$ \cite[Section~$8.4$]{bhatt_absolute_2022}, as constructed in \cite[Construction~$5.8$]{bouis_motivic_2024}. Taking fibres along the horizontal maps in the natural commutative diagram
$$\begin{tikzcd}
    \Z/p^k(i)^{\text{syn}}(X) \ar[r] \arrow[d,"\text{id}"] & \big(L_{\text{Nis}} \tau^{>i} \Z/p^k(i)^{\text{syn}}\big)(X) \arrow[d,"L_{\text{cdh}}"] \\
    \Z/p^k(i)^{\text{syn}}(X) \ar[r] & \big(L_{\text{cdh}} \tau^{>i} \Z/p^k(i)^{\text{syn}}\big)(X),
\end{tikzcd}$$
Theorem~\ref{theoremreview}\,$(3)$ then induces, for any qcqs scheme $X$, a natural map
\begin{equation}\label{diag:BLcomparisonmap}
    \big(L_{\text{Nis}} \tau^{\leqslant i} \Z/p^k(i)^{\text{syn}}\big)(X) \longrightarrow \Z/p^k(i)^{\text{mot}}(X)
\end{equation}
in the derived category $\mathcal{D}(\Z/p^k)$, which we call the {\it Beilinson--Lichtenbaum comparison map}. When $X$ is smooth over a Dedekind domain, using the Beilinson--Lichtenbaum conjecture for Bloch's cycle complexes (\cite{geisser_motivic_2004}), this comparison map is the reduction mod $p^k$ of the comparison map from Bloch's cycle complexes to the motivic complexes $\Z(i)^{\text{mot}}(X)$ (\cite[Definition~$3.23$]{bouis_motivic_2024}). More generally, we will use this 
Beilinson--Lichtenbaum comparison map to formulate our Theorem~\ref{theoremBLPrüfer}. Also note that, for any qcqs scheme~$X$, the composite of the Beilinson--Lichtenbaum comparison map with the syntomic realisation map
$$\big(L_{\text{Nis}} \tau^{\leqslant i} \Z/p^k(i)^{\text{syn}}\big)(X) \longrightarrow \Z/p^k(i)^{\text{mot}}(X) \longrightarrow \Z/p^k(i)^{\text{syn}}(X)$$
is an isomorphism in degrees less than or equal to $i$, since both maps of this composite are actually isomorphisms in this range by Theorem~\ref{theoremreview}$\,(3)$. 

A key point in our arguments will be to use results in classical motivic cohomology to understand the behaviour of the motivic complexes $\Z(i)^{\text{mot}}$ at valuation rings. 

\begin{remark}[Motivic cohomology of valuation rings, \ccite{bouis_motivic_2026}{Section~$5$}]\label{rem:motivic_cohomology_valuation_rings}\label{remarkmotiviccohomologyvaluationrings}
    For every henselian valuation ring $V$ and every integer $i \geqslant 0$, the motivic complex $\Z(i)^{\text{mot}}(V) \in \mathcal{D}(\Z)$ lives in degrees at most $i$.
\end{remark}

In order to compare the Milnor range of integral motivic cohomology with its mod $p^k$ counterpart, we will need the following vanishing result.

\begin{remark}[Motivic cohomology of henselian local rings, \ccite{bouis_weibel_2025}{Section~$2$}]\label{remarkmotiviccohomologyhenselianlocal}
    For every henselian local ring $R$ and every integer $i \geqslant 0$, the motivic cohomology group $\text{H}^j_{\text{mot}}(R,\Z(i)) \colonequals \text{H}^j(\Z(i)^{\text{mot}}(R))$ is zero for $j=i+1$ (\cite[Corollary~$2.10$]{bouis_weibel_2025}), and torsion for $i < j \leqslant 2i$ (\cite[Proposition~$2.4$]{bouis_weibel_2025}).
\end{remark}

\begin{remark}[\'Etale motivic cohomology]\label{rem:etale_motivic_complex}
    Although algebraic $K$-theory of schemes does {\it not} satisfy étale descent (see for instance \cite{clausen_hyperdescent_2021}), one may consider two of its étale-local variants. The first one is étale $K$-theory $\text{K}^{\text{ét}}$, defined here as the étale sheafification of algebraic $K$-theory. Alternatively, one may also consider Selmer $K$\nobreakdash-theory $\text{K}^{\text{Sel}}$, introduced by Clausen \cite{clausen_K-theoretic_2017} and Clausen--Mathew \cite{clausen_hyperdescent_2021} as the pullback $L_{\text{KU}} \text{K} \times_{L_{\text{KU}} \text{TC}} \text{TC}$, where $L_{\text{KU}}$ denotes the Bousfield localisation at the complex $K$-theory spectrum $\text{KU}$. Selmer $K$-theory, has the advantage of being defined at the level of categories, and is even a localising invariant in the sense of \cite{blumberg_universal_2013}. Moreover, it admits a natural map $\text{K}^{\text{ét}} \rightarrow \text{K}^{\text{Sel}}$ which is an isomorphism on homotopy in degrees $\geqslant -1$ (\cite[Theorem~$1.1$]{clausen_hyperdescent_2021}), and is the étale localisation of algebraic $K$-theory when seen in the category of $\mathbb{P}^1$-motivic spectra (\cite[Theorems~$0.1.1$ and~$5.4.4$]{annala_motivic_2023}). 

    Similarly, at the level of motivic cohomology, the motivic complexes $\Z(i)^{\text{mot}}$ are represented by a $\mathbb{P}^1$\nobreakdash-motivic spectrum $\text{H}\!\Z^{\text{mot}}$ (\cite[Corollary~$4.20$]{bouis_weibel_2025}), and one may consider either the étale motivic complexes $\Z(i)^{\text{ét}}$, defined as the étale sheafification of the motivic complexes $\Z(i)^{\text{mot}}$, or the Selmer motivic complexes $\Z(i)^{\text{Sel}}$, defined as the cohomology represented by the étale localisation $\text{H}\!\Z^{\text{Sel}}$ in $\mathbb{P}^1$-motivic spectra of $\text{H}\!\Z^{\text{mot}}$. By construction, these two cohomology theories correspond to the shifted graded pieces of corresponding filtrations on étale $K$-theory and on Selmer $K$-theory, respectively. The Selmer motivic complexes $\Z(i)^{\text{Sel}}$ admit a natural comparison map $\Z(i)^{\text{ét}} \rightarrow \Z(i)^{\text{Sel}}$, and are naturally identified, with mod $p^k$ coefficients, with Bhatt--Lurie's syntomic complexes $\Z/p^k(i)^{\text{syn}}$ (\cite[Remark~$8.4.3$]{bhatt_absolute_2022} and \cite[Example~$9.21$]{annala_atiyah_2024}).
\end{remark}

\section{Proof of the Beilinson--Lichtenbaum phenomenon}\label{sectionBeilinsonLichtenbaum}

\vspace{-\parindent}
\hspace{\parindent}

In this section, we prove Theorems~\ref{theoremintroBL} and~\ref{theoremintroétale}, using Theorem~\ref{theoremGersteninjectivity}. Let us fix a prime number $p$ throughout the section.

\begin{theorem}[Beilinson--Lichtenbaum over Prüfer rings]\label{theoremBLPrüfer}
    Let $X$ be an ind-smooth scheme over a Prüfer ring $P$. Then for all integers $i \geqslant 0$ and $k \geqslant 1$, the 
    Beilinson--Lichtenbaum comparison map \eqref{diag:BLcomparisonmap}
    $$\big(L_{\emph{Nis}} \tau^{\leqslant i} \Z/p^k(i)^{\emph{syn}}\big)(X) \longrightarrow \Z/p^k(i)^{\emph{mot}}(X)$$
    is an equivalence in the derived category $\mathcal{D}(\Z/p^k)$.
\end{theorem}

\begin{proof}
    The presheaves $L_{\text{Nis}} \tau^{\leqslant i} \Z/p^k(i)^{\text{syn}}$ and $\Z/p^k(i)^{\text{mot}}$ are finitary Nisnevich sheaves, so it suffices to prove the result on henselian local ind-smooth $P$-algebras. For such a ring $R$, the statement amounts to proving that the natural map
    $$\tau^{\leqslant i} \Z/p^k(i)^{\text{syn}}(R) \longrightarrow \Z/p^k(i)^{\text{mot}}(R)$$
    is an equivalence in the derived category $\mathcal{D}(\Z/p^k)$. By Theorem~\ref{theoremreview}\,$(3)$, this map is an equivalence in degrees less than or equal to $i$, so it suffices to prove that for every integer $j > i$, the abelian group $\text{H}^j(\Z/p^k(i)^{\text{mot}}(R))$ is zero. 
    
    Let $j \in \Z$ be an integer such that $j > i$. By Remark \ref{rem:motivic_cohomology_valuation_rings}, this vanishing holds if $R$ is a henselian valuation ring. Moreover, the presheaves $\Z/p^k(i)^{\text{mot}}$ are finitary Nisnevich sheaves that satisfy the $\mathbb{P}^1$-bundle formula (Theorem \ref{theoremreview}\,(1)), and thus, are deflatable (Lemma~\ref{lemmapbfimpliesdeflatable}). The vanishing in general is then a consequence of Gersten's injectivity (Theorem~\ref{theoremGersteninjectivity}).
\end{proof}

\begin{remark}\label{remarkBLfalseingeneral}
    The statement of Theorem~\ref{theoremBLPrüfer} is not expected to hold over an arbitrary base ring. Indeed, such a statement would imply, via the Atiyah--Hirzebruch spectral sequence, that non-connective algebraic $K$-theory (with finite coefficients) is the Nisnevich sheafification of connective algebraic $K$-theory on arbitrary schemes, which is known to be true only after also enforcing the projective bundle formula \cite{annala_motivic_2023}. However, we expect that the statement of Theorem~\ref{theoremBLPrüfer} should also hold over an arbitrary regular base ring. 
\end{remark}

\begin{corollary}\label{corollarymotivicoverPrüferlowdegrees}
    Let $X$ be an ind-smooth scheme over a Prüfer ring. Then for all integers $i \geqslant 0$ and $k \geqslant 1$, the motivic complex $\Z/p^k(i)^{\emph{mot}}(X)$ lives, Nisnevich-locally on $X$, in degrees at most $i$.
\end{corollary}

\begin{proof}
    This is a consequence of Theorem~\ref{theoremBLPrüfer}.
\end{proof}

\begin{corollary}\label{propositionP1bundleformulaforNislissefinitecoefficients}
    For every Prüfer ring $P$ and every integer $k \geqslant 1$, the family of presheaves $$\{L_{\emph{Nis}} \tau^{\leqslant i} \Z/p^k(i)^{\emph{syn}}\}_{i \geqslant 0}$$ satisfies the projective bundle formula on ind-smooth $P$-schemes.
\end{corollary}

\begin{proof}
    The family of presheaves $\{\Z/p^k(i)^{\text{mot}}\}_{i \geqslant 0}$ satisfies the projective bundle formula on all qcqs schemes (Theorem \ref{theoremreview}\,(1)), so this is a consequence of Theorem~\ref{theoremBLPrüfer}.
\end{proof}

\begin{corollary}
    Let $f : P \rightarrow P'$ be a morphism of Prüfer rings, and $\emph{H}\!\Z^{\emph{mot}}_P \in \emph{MS}(P)$ and $\emph{H}\!\Z^{\emph{mot}}_{P'} \in \emph{MS}(P')$ be the $\mathbb{P}^1$-motivic spectra representing the motivic complexes $\Z(i)^{\emph{mot}}$ over $P$ and $P'$ respectively. Then the natural map
    $$f^\ast \emph{H}\!\Z^{\emph{mot}}_P \longrightarrow \emph{H}\!\Z^{\emph{mot}}_{P'}$$
    is an equivalence in $\emph{MS}(P')$.
\end{corollary}

\begin{proof}
    It suffices to prove the result rationally, and modulo $p$ for every prime number $p$. The result modulo~$p$ is a consequence of Theorem~\ref{theoremBLPrüfer} and Corollary~\ref{propositionP1bundleformulaforNislissefinitecoefficients}, where we use the fact that the functor $\tau^{\leqslant i} \F_p(i)^{\text{syn}} : \text{Rings} \rightarrow \mathcal{D}(\Z/p^k)$ is left Kan extended from smooth $\Z$-algebras (\cite[Proposition~$2.7$]{bouis_weibel_2025}). The result rationally holds true for any morphism of commutative rings $f : P \rightarrow P'$, as a consequence of the Adams decomposition (Theorem~\ref{theoremreview}\,$(2)$) and of base change for algebraic $K$-theory in $\mathbb{P}^1$-motivic spectra (\cite[Theorem~$0.1.1$]{annala_motivic_2023}).
\end{proof}

In the next results, we consider the étale motivic complexes $\Z(i)^{\text{ét}}$ and $\Z(i)^{\text{Sel}}$ introduced in Remark~\ref{rem:etale_motivic_complex}.

\begin{corollary}\label{corollaryétale=syn}
    Let $X$ be an ind-smooth scheme over a Prüfer ring. Then for all integers $i \geqslant 0$ and $k \geqslant 1$, the natural map
    $$\Z/p^k(i)^{\emph{ét}}(X) \longrightarrow \Z/p^k(i)^{\emph{syn}}(X)$$
    is an equivalence in the derived category $\mathcal{D}(\Z/p^k)$.
\end{corollary}

\begin{proof}
    By \cite[Theorem~$5.1\,(1)$ and Corollary~$5.43$]{antieau_beilinson_2022}, the weight-$i$ Bhatt--Morrow--Scholze's syntomic complex lives, étale-locally on any qcqs scheme $X$, in degrees at most $i$. By \cite[Remark~$8.4.4$]{bhatt_absolute_2022}, the syntomic complex $\Z/p^k(i)^{\text{syn}}$ then satisfies the same property, hence the desired result is a consequence of Theorem~\ref{theoremBLPrüfer}.
\end{proof}

\begin{proposition}
    Let $X$ be an ind-smooth scheme over a Prüfer ring. Then for every integer $i \geqslant 0$, the natural map
    $$\Z(i)^{\emph{ét}}(X) \longrightarrow \Z(i)^{\emph{Sel}}(X)$$
    is an equivalence in the derived category $\mathcal{D}(\Z)$.
\end{proposition}

\begin{proof}
    It suffices to prove the result rationally, and modulo $p$ for every prime number $p$. The result modulo~$p$ is Corollary~\ref{corollaryétale=syn}. The result rationally holds true for any qcqs scheme $X$, as a consequence of the Adams decomposition (Theorem~\ref{theoremreview}\,$(2)$) and of the fact that the natural map $\text{K}(X) \rightarrow \text{K}^{\text{Sel}}(X)$ is a rational equivalence (\cite[Example~$6.2$]{clausen_hyperdescent_2021}).
\end{proof}

\begin{theorem}
    Let $X$ be an ind-smooth scheme over a Prüfer ring. Then for every integer $i \geqslant 0$, the motivic complex $\Z(i)^{\emph{mot}}(X)$ lives, Nisnevich-locally on $X$, in degrees at most $i$, and the induced natural map
    $$\Z(i)^{\emph{mot}}(X) \longrightarrow \big(L_{\emph{Nis}} \tau^{\leqslant i} \Z(i)^{\emph{ét}}\big)(X)$$ 
    is an equivalence in the derived category $\mathcal{D}(\Z)$.
\end{theorem}

\begin{proof} 
    Let $R$ be the henselisation of a local ring of $X$. The ring $R$ is then a henselian local ind-smooth algebra over a valuation ring, so the motivic complex $\Z/p^k(i)^{\text{mot}}(R)$ lives in degrees at most $i$ (Corollary~\ref{corollarymotivicoverPrüferlowdegrees}). By the short exact sequence of abelian groups
    $$0 \longrightarrow \text{H}^j_{\text{mot}}(R,\Z(i))/p^k \longrightarrow \text{H}^j_{\text{mot}}(R,\Z/p^k(i)) \longrightarrow \text{H}^{j+1}_{\text{mot}}(R,\Z(i))[p^k] \longrightarrow 0,$$
    for every prime number $p$ and every integer $j \geqslant 0$, this implies that the abelian group $\text{H}^j_{\text{mot}}(R,\Z(i))$ is torsionfree for $j \geqslant i+2$. By Remark~\ref{remarkmotiviccohomologyhenselianlocal}, where we use that $R$ is henselian local, the abelian group $\text{H}^{i+1}_{\text{mot}}(R,\Z(i))$ is zero, hence torsionfree. It then suffices to prove that the rational motivic complex $\Q(i)^{\text{mot}}(R)$ lives in degrees at most $i$. This complex vanishes in degrees $i+1$ to $2i$ by Remark~\ref{remarkmotiviccohomologyhenselianlocal}, and is a direct summand of the shift $\text{K}(R)_{\Q}[-2i]$ of the rational algebraic $K$\nobreakdash-theory of~$R$ (Theorem~\ref{theoremreview}\,$(2)$). Moreover, by \cite[Corollary~$2.3$ and Proposition~$2.4.1$]{antieau_K-theory_2022}, the negative $K$\nobreakdash-groups of ind-smooth algebras over valuation rings are zero. In particular, the negative $K$\nobreakdash-groups of $R$ are zero, and the cohomology groups of the complex $\Q(i)^{\text{mot}}(R)$ also vanish in degrees more than $2i$. This proves the first claim. In particular, the canonical map $\Z(i)^{\text{mot}}(X) \rightarrow \Z(i)^{\text{ét}}(X)$ naturally factors as
    $$\Z(i)^{\text{mot}}(X) \longrightarrow \big(L_{\text{Nis}} \tau^{\leqslant i} \Z(i)^{\text{ét}}\big)(X) \longrightarrow \Z(i)^{\text{ét}}(X)$$
    in the derived category $\mathcal{D}(\Z)$. 
    
    We prove now that the first map of this composite is an equivalence. It suffices to do so rationally and modulo $p$ for every prime number $p$. By Thomason (\cite[Proposition~$2.14$]{thomason_algebraic_1985}), rational algebraic $K$-theory satisfies étale descent, therefore, the same holds for rational motivic cohomology (Theorem~\ref{theoremreview}\,$(2)$). The claim that the natural map $\Q(i)^{\text{mot}}(X) \rightarrow \big(L_{\text{Nis}} \tau^{\leqslant i} \Q(i)^{\text{mot}}\big)(X)$ is an equivalence is then a consequence of the previous paragraph. Modulo $p$, and using again that the motivic complex $\Z(i)^{\text{mot}}$ vanishes in \hbox{degree $i+1$} for henselian local rings (Remark~\ref{remarkmotiviccohomologyhenselianlocal}), the desired statement is equivalent to the fact that the natural map
    $$\F_p(i)^{\text{mot}}(X) \longrightarrow \big(L_{\text{Nis}} \tau^{\leqslant i} \F_p(i)^{\text{ét}}\big)(X)$$
    is an equivalence in the derived category $\mathcal{D}(\F_p)$. After using Corollary~\ref{corollaryétale=syn} to express the $\F_p(i)^{\text{ét}}$ in terms of syntomic cohomology, this is a consequence of Theorem~\ref{theoremBLPrüfer}.
\end{proof}

\begin{proof}[Proof of Theorem~\ref{theoremintroétale}]
    The henselian local ring $R$ is local for the Nisnevich topology, and is ind-smooth over a valuation ring $V$ by definition, so the Beilinson--Lichtenbaum comparison map induces a natural equivalence
    $$\Z/p^k(i)^{\text{mot}}(R) \simeq \tau^{\leqslant i} \Z/p^k(i)^{\text{syn}}(R)$$
    in the derived category $\mathcal{D}(\Z/p^k)$ (Theorem~\ref{theoremBLPrüfer}). If the commutative ring $V$ is $F$-smooth, then ind-smooth algebras over $V$ are $F$-smooth (\cite[Propositions~$4.6$ and~$4.9$]{bhatt_syntomic_2023}). The desired identification of the complex $\tau^{\leqslant i} \Z/p^k(i)^{\text{syn}}(R)$ in terms of étale cohomology, for $p$-torsionfree $F$-smooth rings $R$, is then \cite[Theorem~$1.8$]{bhatt_syntomic_2023}, whose proof relies on absolute prismatic cohomology.
\end{proof}

We end this section with two consequences of Theorem~\ref{theoremBLPrüfer}, which will be used in Remark~\ref{remarkcrystallineconjecture}.

\begin{corollary}\label{corollarycomparisongenericfibre}
    Let $V$ be a $p$-adically separated valuation ring of residue characteristic $p$ whose $p$-completion is perfectoid ({\it e.g.}, a perfectoid valuation ring of residue characteristic $p$), $F$ be the fraction field of $V$, and $X$ be an ind-smooth scheme over $V$. Then for all integers $i \geqslant 0$ and $k \geqslant 1$, the natural map
    $$\Z/p^k(i)^{\emph{mot}}(X) \longrightarrow \Z/p^k(i)^{\emph{mot}}(X_F)$$
    is an equivalence in the derived category $\mathcal{D}(\Z/p^k)$.
\end{corollary}

\begin{proof}
    If $V$ is an $\F_p$-algebra, then $V$ is a filtered colimit of smooth $\F_p$-algebras, and the result is a consequence of Geisser--Levine's description of the $p$-adic motivic cohomology of smooth $\F_p$-schemes (see \cite[Remark~$3.5$]{antieau_K-theory_2022}). We assume now that $p$ is nonzero in $V$. In particular, because $V$ is $p$-adically separated, we know that $V[\tfrac{1}{p}]=F$ (\ccite{fujiwara_kato_foundations_rigid_geometry}{Chapter~0, Proposition~6.7.2}), and it then suffices to prove that the natural map
    $$\Z/p^k(i)^{\text{mot}}(X) \longrightarrow \Z/p^k(i)^{\text{mot}}(X[\tfrac{1}{p}])$$
    is an equivalence in the derived category $\mathcal{D}(\Z/p^k)$. Perfectoid rings are $F$-smooth (\cite[Proposition~$4.12$]{bhatt_absolute_2022}), so $V$ is a $p$-torsionfree $F$-smooth valuation ring. By Theorem~\ref{theoremintroétale}, it then suffices to prove that the natural map $\text{H}^i_{\text{mot}}(X,\Z/p^k(i)) \rightarrow \text{H}^i_{\text{mot}}(X[\tfrac{1}{p}],\Z/p^k(i))$ is an isomorphism of abelian groups. The presheaves $\Z/p^k(i)^{\text{mot}}(-)$ and $\Z/p^k(i)^{\text{mot}}(-[\tfrac{1}{p}])$ are deflatable finitary Nisnevich sheaves on qcqs $V$-schemes, so it suffices to prove that the natural map $$\text{H}^i_{\text{mot}}(V',\Z/p^k(i)) \longrightarrow \text{H}^i_{\text{mot}}(V'[\tfrac{1}{p}],\Z/p^k(i))$$ is an isomorphism for every henselian valuation ring $V'$ which is ind-smooth over $V$ (Theorem~\ref{theoremGersteninjectivity}). This map of abelian groups is injective by Theorem~\ref{theoremreview}\,$(3)$ and \cite[Corollary~$3.2$]{bouis_cartier_2023} (whose proof holds for general valuation rings extensions $V'$ of a perfectoid valuation ring $V$, but is considerably simplified when $V'$ is ind-smooth over~$V$). The symbol maps in motivic cohomology (\cite[Section~$2.2$]{bouis_weibel_2025}) naturally fit in a commutative diagram
    $$\begin{tikzcd}
        (V'^{\times}/p^k)^{\otimes i} \ar[r] \ar[d] & \text{H}^i_{\text{mot}}(V',\Z/p^k(i)) \ar[d] \\
        (V'[\tfrac{1}{p}]^{\times}/p^k)^{\otimes i} \ar[r] & \text{H}^i_{\text{mot}}(V'[\tfrac{1}{p}],\Z/p^k(i))
    \end{tikzcd}$$
    of abelian groups. The bottom horizontal map is surjective by \cite[Theorem~$14.1$]{bloch_p-adic_1986}. The value group $\Gamma_{V'} \colonequals \text{coker}(V'^{\times} \rightarrow V'[\tfrac{1}{p}]^{\times})$ of the valuation ring $V'$ is naturally identified with that of the valuation ring~$V$ (see \cite[Lemma~$3.10$]{phd-thesis}). The value group of the $p$-adic completion of $V$ is naturally identified with that of its tilt $V^{\flat}$, which is $p$-divisible, hence so is the value group of $V$ (\ccite{fujiwara_kato_foundations_rigid_geometry}{Chapter~0, Theorem~9.1.1}). This implies that the left vertical map is an isomorphism for $i=1$, and hence for all $i \geqslant 0$. In particular, the right vertical map is surjective, which concludes the proof.
\end{proof}

\begin{remark}\label{rem:annala_elmanto}
    Using Levine's absolute purity for motivic cohomology over discrete valuation rings, Corollary~\ref{corollarycomparisongenericfibre} was also recently proved independently by Annala and Elmanto in the special cases where $V$ is either $\Z_p[p^{1/p^{\infty}}]$ or $\Z[\zeta_{p^{\infty}}]$ \cite[Theorem~$3.1$ and Variant~$3.4$]{annala_elmanto}.
\end{remark}

For smooth schemes over an algebraically closed field in which $p$ is invertible, the following result was first proved by Suslin. We denote by $p$-$\text{dim}(R)$ the $p$-dimension of a $\F_p$-algebra $R$, and refer to \cite[Section~$2$]{kerz_towards_2021} for the associated results needed in the following proof. In particular, note that if $\kappa$ is a perfect $\F_p$-algebra, then $p\text{-dim}(\kappa)=0$.

\begin{corollary}\label{propositionmotivicsyntomicinhighweights}
	Let $X \rightarrow S$ be an affine smooth morphism of schemes of relative dimension $d$. Then for every integer $k \geqslant 1$, the syntomic realisation map
	$$\Z/p^k(i)^{\emph{mot}}(X) \longrightarrow \Z/p^k(i)^{\emph{syn}}(X)$$
	is an equivalence in the derived category $\mathcal{D}(\Z/p^k)$ in each of the following cases:
	\begin{enumerate}
		\item $S$ is the spectrum of a separably closed field $K$ of characteristic different from $p$, and $i \geqslant d$;
		\item $S$ is the spectrum of a field $\kappa$ of characteristic $p$, and $i \geqslant d + p\emph{-dim}(\kappa) + 1$;
		\item $S$ is the spectrum of a rank one valuation ring $V$ of mixed characteristic $(0,p)$ with separably closed fraction field $K$ and residue field $\kappa$, and $i \geqslant d + p\emph{-dim}(\kappa) + 1$.
	\end{enumerate}
\end{corollary}

\begin{proof}
    By Theorem~\ref{theoremBLPrüfer}, it suffices to prove in each case that the syntomic complex $\Z/p^k(i)^{\text{syn}}$ lives, Nisnevich-locally on $X$, in degrees at most $i$.
    
    $(1)$ If $p$ is invertible in the separably closed field $K$, then this is a consequence of the Artin--Grothendieck vanishing theorem \cite[exposé~XIV, corollaire~$3.2$]{SGA4iii}, which states in particular that the étale cohomology $R\Gamma_{\text{ét}}(R,\mu_{p^k}^{\otimes i})$ of any finite type $K$-algebra $R$ of dimension $d$ lives in degrees at most $d$, hence, in particular, in degrees at most $i$.

    $(2)$ If $\kappa$ is of characteristic $p$, then $\tau^{>i} \Z/p^k(i)^{\text{syn}}(R) \simeq \widetilde{\nu}_k(i)(R)[-i-1]$ for every $\kappa$-algebra $R$ (\cite[Theorem~$5.1$ and Corollary~$5.43$]{antieau_beilinson_2022}). If $R$ is moreover of finite type over $\kappa$, then the $p$-dimension of $R$ is equal to $d+p\text{-dim}(\kappa)$. This implies that $\Omega^i_{R/\kappa}$ vanishes for every integer $i \geqslant d + p\text{-dim}(\kappa) + 1$, hence so does $\widetilde{\nu}_k(i)(R)$.

    $(3)$ Let $R$ be the henselisation of a local ring of $X$. If $p$ is invertible in $R$, then $R$ is the henselisation of a local ring of $X[\tfrac{1}{p}]$, which is a smooth scheme of dimension $d$ over the fraction field $K$ of $V$, hence the claim reduces to $(1)$. If $p$ is not invertible in the henselian local ring $R$, then in particular the ring $R$ is $p$-henselian, and its syntomic cohomology is naturally identified with Bhatt--Morrow--Scholze's syntomic cohomology (\cite[Remark~$8.4.4$]{bhatt_absolute_2022}). In particular, the natural map 
    $$\tau^{>i} \Z/p^k(i)^{\text{syn}}(R) \longrightarrow \tau^{>i} \Z/p^k(i)^{\text{syn}}(R \otimes_V \kappa)$$
    is an equivalence in the derived category $\mathcal{D}(\Z/p^k)$ (\cite[Theorem~$5.2$]{antieau_beilinson_2022}, see \cite[Theorem~$2.27$]{bouis_motivic_2024} for the precise version that we use). The ring $R/\mathfrak{m}R$ is the henselisation of a local ring of $X_{\kappa}$, which is a scheme of relative dimension $d$ over the field $\kappa$, hence the claim again reduces to $(2)$.    
\end{proof}

\section{$\mathbb{A}^1$-invariance and comparison with Bloch's cycle complexes}\label{sectionA1invariance}

\vspace{-\parindent}
\hspace{\parindent}

In this section, we prove that the motivic complexes $\Z(i)^{\text{mot}}(X)$ recover the classical motivic complexes $$\Z(i)^{\text{cla}}(X) \colonequals \big(L_{\text{Zar}}\, z^i(-,\bullet)\big)(X)[-2i]$$ on smooth schemes $X$ over a Dedekind domain (Theorem~\ref{theoremclassicalmotiviccomparison}). To formulate this result, we use the {\it classical-motivic comparison map} $\Z(i)^{\text{cla}}(X) \rightarrow \Z(i)^{\text{mot}}(X)$ of \cite[Definition~$3.23$]{bouis_motivic_2024}. We also denote by $\Z(i)^{\text{cla}}$ the unique functorial extension of the previous Bloch cycle complex $\Z(i)^{\text{cla}}$ to schemes that are ind-smooth over a Dedekind domain.

\begin{theorem}[Comparison to classical motivic cohomology]\label{theoremclassicalmotiviccomparison}
Let $X$ be an ind-smooth scheme over a Dedekind domain $B$. Then for every integer $i \geqslant 0$, the classical-motivic comparison map
$$\Z(i)^{\emph{cla}}(X) \longrightarrow \Z(i)^{\emph{mot}}(X)$$
is an equivalence in the derived category $\mathcal{D}(\Z)$.
\end{theorem}

\begin{proof}
    If the Dedekind domain $B$ contains a field, then it is a filtered colimit of smooth algebras over this field, and the result is a consequence of \cite[Corollary~$6.4$]{elmanto_motivic_2023}. We assume now that the Dedekind domain $B$ is of mixed characteristic. It suffices to prove the result rationally, and modulo $p$ for every prime number~$p$.

    The result rationally is a consequence of the rational splitting of algebraic $K$-theory induced by Adams operations. More precisely, we use the splitting induced by \cite[Lemma~4.11]{bouis_motivic_2024} for the filtrations $\text{Fil}^\star_{\text{cla}} \text{K}(-; \Q)$ (which is $\N$-indexed by construction) and $\text{Fil}^\star_{\text{mot}} \text{K}(-; \Q)$ (which is $\N$-indexed by \cite[Proposition~4.46]{bouis_motivic_2024}). These decompositions are compatible with the classical-motivic comparison map because of the compatibility between the associated Adams operations (\cite[Section~4.1]{bouis_motivic_2024}).
    
    Let $p$ be a prime number. By \cite[Corollary~$4.4$]{geisser_motivic_2004}, Nisnevich-locally, the classical motivic complex $\F_p(i)^{\text{cla}}$ lives in degrees at most $i$, so the natural composite map
    $$\F_p(i)^{\text{cla}}(X) \longrightarrow \F_p(i)^{\text{mot}}(X) \longrightarrow \F_p(i)^{\text{syn}}(X)$$
    naturally factors through a composite map
    $$\F_p(i)^{\text{cla}}(X) \longrightarrow \big(L_{\text{Nis}} \tau^{\leqslant i} \F_p(i)^{\text{mot}}\big)(X) \longrightarrow \big(L_{\text{Nis}} \tau^{\leqslant i} \F_p(i)^{\text{syn}}\big)(X)$$
    in the derived category $\mathcal{D}(\F_p)$. By \cite[Theorems~$1.2\,(2)$ and~$1.3$]{geisser_motivic_2004} and \cite[Theorem~$5.8$ and Example~$5.9$]{bhatt_syntomic_2023}, this composite map is an equivalence. Moreover, by Corollary~\ref{corollarymotivicoverPrüferlowdegrees}, the middle term is naturally identified with the motivic complex $\F_p(i)^{\text{mot}}(X)$. The second map is then an equivalence by Theorem~\ref{theoremBLPrüfer}, so the natural map 
    $$\F_p(i)^{\text{cla}}(X) \longrightarrow \F_p(i)^{\text{mot}}(X)$$
    is also an equivalence in the derived category $\mathcal{D}(\F_p)$, as desired.
\end{proof}
	
\begin{remark}
    When $B$ is a field, Theorem~\ref{theoremclassicalmotiviccomparison} is \cite[Corollary~$6.4$]{elmanto_motivic_2023}.
\end{remark}

\begin{remark}[Fontaine's crystalline conjecture]\label{remarkcrystallineconjecture}
    Fontaine's crystalline conjecture \cite[conjecture~A.$11$]{fontaine_certains_1982}, first proved in general by Faltings \cite{faltings_crystalline_1989}, asserts that for $X$ a smooth proper scheme over a complete discrete valuation ring $\mathcal{O}_K$ of mixed characteristic $(0,p)$ and with perfect residue field $k$, there are natural isomorphisms of abelian groups
    $$\text{H}^n_{\text{ét}}(X_{\overline{K}},\Z_p) \otimes_{\Z_p} B_{\text{crys}} \cong \text{H}^n_{\text{crys}}(X_k/W(k)) \otimes_{W(k)} B_{\text{crys}} \quad \quad  \text{ for all } n \geqslant 0,$$
    compatible with filtrations, Galois actions, and Frobenius operators. As observed by Nizio\l{} \cite{niziol_p-adic_2006}, this comparison should be reminiscent of a suitable theory of motivic cohomology for mixed characteristic schemes. In this remark, we explain how the results of the present paper can be used to realise this expectation.

    By \cite[Section~$4$]{niziol_crystalline_1998} (see also \cite{fontaine_p-adic_1987}, or \cite[$2.4$]{faltings_p-adic_1988}), the key point is to construct, at least for $i \gg 0$, a natural Galois-equivariant equivalence
    $$R\Gamma_{\text{ét}}(X_{\overline{K}},\mu_{p^k}^{\otimes i}) \xlongrightarrow{\sim} \Z/p^k(i)^{\text{syn}}(X_{\mathcal{O}_{\overline{K}}}),$$
    compatible with Poincaré duality and cycle classes, for every integer $k \geqslant 1$.\footnote{In \cite[Theorem~$4.1$]{niziol_crystalline_1998}, and using algebraic $K$-theory, Nizio\l{} constructs such an equivalence for $X$ smooth projective of pure dimension $d$ over $\mathcal{O}_K$, and with the restrictions that $i \geqslant \tfrac{16}{3}d^3 + 8d^2 + \tfrac{14}{3}d$ and $p \geqslant 2i + \tfrac{16}{3}d^3 + 8d^2 + \tfrac{8}{3}d + 2$. While the latter restriction is partly explained by the use of Fontaine--Messing's syntomic cohomology, which works well integrally only in weights $i$ less than $p$, the former restriction comes from Thomason's comparison between algebraic $K$-theory and étale cohomology.} This equivalence can now be seen as the inverse of the bottom map in the commutative diagram
    $$\begin{tikzcd}
        \Z/p^k(i)^{\text{mot}}(X_{\overline{K}}) \ar[d] & \Z/p^k(i)^{\text{mot}}(X_{\mathcal{O}_{\overline{K}}}) \ar[l] \ar[d] \\
        R\Gamma_{\text{ét}}(X_{\overline{K}},\mu_{p^k}^{\otimes i}) & \Z/p^k(i)^{\text{syn}}(X_{\mathcal{O}_{\overline{K}}}) \ar[l]
    \end{tikzcd}$$
    where the horizontal maps are induced by functoriality of the motivic and syntomic complexes, and the vertical maps are instances of the syntomic realisation map \cite[Construction~$5.8$]{bouis_motivic_2024}. More precisely, for every smooth scheme $X_{\mathcal{O}_{\overline{K}}}$ over the valuation ring $\mathcal{O}_{\overline{K}}$, the vertical maps of this diagram are equivalences for $i \geqslant d+1$ by Proposition~\ref{propositionmotivicsyntomicinhighweights}, and the top horizontal map is an equivalence for all $i \geqslant 0$ by Corollary~\ref{corollarycomparisongenericfibre}. The compatibility with cycle classes is then a consequence of Theorem~\ref{theoremclassicalmotiviccomparison} (where we use that $\mathcal{O}_{\overline{K}}$ is a filtered union of discrete valuation rings), and the compatibility with Poincaré duality is in turn a consequence of the compatibility with cycle classes (see \cite[proof of Lemma~$4.2$]{niziol_crystalline_1998}).
\end{remark}

For the rest of this section, we use Theorem~\ref{theoremGersteninjectivity} to study the $\mathbb{A}^1$-invariance of the motivic complexes~$\Z(i)^{\text{mot}}$.

\begin{definition}\label{definitioncdhlocallyFsmooth}
    A commutative ring $R$ is {\it cdh-locally $F$-smooth} if every valuation ring $V$ over $R$ is $F$\nobreakdash-smooth\footnote{$F$-smoothness implicitly depends on a prime number $p$. Here we say that a valuation ring $V$ is $F$-smooth if it is $F$-smooth for every prime number $p$. Note that valuation rings $V$ are local, so that $F$-smoothness is vacuously true for almost all prime numbers $p$, {\it i.e.}, for those that are invertible in $V$.} in the sense of \cite[Definition~$1.7$]{bhatt_syntomic_2023}.
\end{definition}

\begin{examples}
    \begin{enumerate}
        \item A field is cdh-locally $F$-smooth. Indeed, valuation rings of characteristic zero are vacuously $F$-smooth, and valuation rings of positive characteristic are Cartier smooth by results of Gabber--Ramero and Gabber (see \cite[Section~$2$]{kelly_k-theory_2021}), hence $F$-smooth (\cite[Proposition~$4.14$]{bhatt_syntomic_2023}).
        \item A perfectoid valuation ring is cdh-locally $F$-smooth. Indeed, valuation rings over a perfectoid valuation ring are $F$-smooth (\cite{bouis_cartier_2023}).
        \item Any algebra over a cdh-locally $F$-smooth ring is cdh-locally $F$-smooth.
    \end{enumerate}
\end{examples}

\begin{remark}\label{remarkconjecturevaluationringsFsmooth}
    Conjecturally, all valuation rings are $F$-smooth, {\it i.e.}, $\Z$ is cdh-locally $F$-smooth. This conjecture is a consequence of Zariski's local uniformisation conjecture, which states that every valuation ring is a filtered colimit of regular local rings. Indeed, regular rings are $F$-smooth, and the class of $F$-smooth rings is stable under filtered colimits (\cite{bhatt_syntomic_2023}).
\end{remark}

The motivation for Definition~\ref{definitioncdhlocallyFsmooth} is the following result of Bachmann--Elmanto--Morrow, where the cdh-local motivic complex $\Z(i)^{\text{cdh}}$ is the cdh sheafification of the motivic complex $\Z(i)^{\text{mot}}$.

\begin{theorem}[\cite{bachmann_A^1-invariant_2025}]\label{theoremBEM}
    On qcqs schemes over a cdh-locally $F$-smooth ring $B$, the cdh-local motivic complexes $\Z(i)^{\emph{cdh}}$, $i \geqslant 0$, are $\mathbb{A}^1$-invariant and satisfy the $\mathbb{P}^1$-bundle formula.
\end{theorem}

\begin{proof}
    Given a qcqs scheme $X$, Bachmann--Elmanto--Morrow \cite[Section~8.1]{bachmann_A^1-invariant_2025} define an hypothesis called $\text{Val}(X,p,i)$, for $p$ a prime number and $i \ge 0$ an integer. This hypothesis, for every prime number~$p$ and integer $i \ge 0$, is satisfied if $X$ is cdh-locally $F$-smooth (\cite[Theorem~1.8]{bhatt_syntomic_2023}), and implies that the cdh-local motivic complexes $\Z(i)^{\text{cdh}}$ agree with the $\mathbb{A}^1$-invariant motivic complexes $\Z(i)^{\mathbb{A}}$ (\cite[Theorem~1.10]{bachmann_A^1-invariant_2025}), which are known to satisfy the $\mathbb{P}^1$-bundle formula (\cite[Theorem~1.1.7]{bachmann_A^1-invariant_2025}).
\end{proof}

\begin{remark}
    At the level of $K$-theory, the analogue of Theorem~\ref{theoremBEM} holds for arbitrary qcqs schemes. More precisely, by a result of Kerz--Strunk--Tamme \cite{kerz_algebraic_2018}, $KH$-theory is the cdh-sheafification of $K$-theory on qcqs schemes.
\end{remark}

The following result is a motivic analogue of the results of Antieau--Mathew--Morrow \cite[Section~$2$]{antieau_K-theory_2022} on the $\mathbb{A}^1$-invariance of algebraic $K$-theory.

\begin{theorem}\label{theoremA1invariancemain}
    Let $B$ be either a Dedekind domain or a cdh-locally $F$-smooth Prüfer domain, and $X$ be an ind-smooth scheme over $B$. Then for all integers $i \geqslant 0$ and $m \geqslant 1$, the natural map
    $$\Z(i)^{\emph{mot}}(X) \longrightarrow \Z(i)^{\emph{mot}}(\mathbb{A}^m_X)$$
    is an equivalence in the derived category $\mathcal{D}(\Z)$.
\end{theorem}

\begin{proof}
     If $B$ is a Dedekind domain, this is a consequence of Theorem~\ref{theoremclassicalmotiviccomparison} and of the $\mathbb{A}^1$-invariance of Bloch cycle complexes. 
     
     Assuming now that $B$ is a cdh-locally $F$-smooth Prüfer domain. In this case, it suffices to prove that the canonical map $\Z(i)^{\text{mot}} \rightarrow \Z(i)^{\text{cdh}}$ is an equivalence on ind-smooth schemes over $B$ (Theorem~\ref{theoremBEM}). Since the Nisnevich finitary sheaves $\Z(i)^{\text{mot}}$ and $\Z(i)^{\text{cdh}}$ are both deflatable (Theorems~\ref{theoremreview}\,$(1)$ and~\ref{theoremBEM}), the same is true for the fibre $$\scr{F}(-)\colonequals \text{fib}\big(\Z(i)^{\text{mot}}(-) \rightarrow \Z(i)^{\text{cdh}}(-)\big).$$ 
     Hence, by Corollary~\ref{cor:theoremGersteninjectivity}, it suffices to show that $\scr{F}(V)=0$ for any henselian valuation ring $V$ over $B$. To prove this, note that henselian valuation rings are the local rings for the cdh topology, and that the presheaf $\Z(i)^{\text{cdh}}$ is the cdh sheafification of the presheaf $\Z(i)^{\text{mot}}$. In particular, for every henselian valuation ring $V$, the natural map
     $$\Z(i)^{\text{mot}}(V) \longrightarrow \Z(i)^{\text{cdh}}(V)$$
     is an equivalence in the derived category $\mathcal{D}(\Z)$, {\it i.e.}, $\mathscr{F}(V)=0$.
\end{proof}
	
\bibliographystyle{alpha}
	
{\footnotesize
\bibliography{biblio.bib}
}

\medskip

{\footnotesize 
\textsc{School of Mathematics, Institute for Advanced Study, 08540 Princeton, New Jersey, USA} \par  
  \textit{Email address:} \texttt{tbouis@ias.edu} \par
  \textit{URL:} \texttt{https://tessbouis.com} \par
  
  \addvspace{\medskipamount}

  \textsc{St.~George Campus, University of Toronto, Toronto ON M5S 2E4, Canada} \par  
  \textit{Email address:} \texttt{arnab.kundu@utoronto.ca} \par
  \textit{URL:} \texttt{https://www.arnabkundu.com} \par
}

\end{document}